\documentclass[a4paper,fleqn]{cas-dc}

\usepackage[authoryear,longnamesfirst]{natbib}
\usepackage{booktabs}
\usepackage{caption}
\usepackage{graphicx,subfig}
\usepackage{amsmath,lipsum}
\usepackage{cuted}
\usepackage{graphicx}
\usepackage{subfig}
\usepackage{subfloat}
\usepackage{epstopdf}
\usepackage{algorithm}
\usepackage{algorithmic}
\usepackage{caption}

\newtheorem{remark}{Remark}[section]

\def\tsc#1{\csdef{#1}{\textsc{\lowercase{#1}}\xspace}}
\tsc{WGM}
\tsc{QE}
\tsc{EP}
\tsc{PMS}
\tsc{BEC}
\tsc{DE}

\begin{document}
\let\WriteBookmarks\relax
\def\floatpagepagefraction{1}
\def\textpagefraction{.001}

\shorttitle{A novel model-based parameters estimation combining local optimization and global optimization of nonlinear ship models with physical experiment dataset}


\title [mode = title]{A novel model-based parameters estimation combining local optimization and global optimization of nonlinear ship models with physical experiment dataset}                      
\tnotemark[1]

\tnotetext[1]{This research was funded by the National Key R\&D Program of China (2022YFE0125200), National Natural Science Foundation of China (62003250, 52272425), and Chinese Scholarship Council (CSC NO. 202206950043).}

\author[1,2]{Xu You}[style = Chinese, orcid=0000-0002-4141-8222]
\credit{Conceptualization of this study, Methodology, Software, Writing - Original Draft}

\address[1]{State Key Laboratory of Maritime Technology and Safety, Wuhan University of Technology, Wuhan, 430063, China}
\address[2]{School of Transportation and Logistics Engineering, Wuhan University of Technology, Wuhan, China}
\address[3]{Intelligent Transportation Systems Research Center, Wuhan University of Technology, Wuhan, China}
\address[4]{National Engineering Research Center for Water Transport Safety, Wuhan University of Technology, Wuhan, China}
\address[5]{Department of Mechanical Engineering, University College London, Torrington Place, London, UK.}

\author[1,2,3,4]{Xinping Yan}[style = Chinese, orcid=0000-0002-2265-2689]
\credit{Writing - Review \& Editing, Supervision, Funding acquisition}

\author[1,3,4]{Jialun Liu}[style=Chinese, orcid=0000-0002-0802-8986]
\ead{jialunliu@whut.edu.cn}
\credit{Writing - Review \& Editing, Supervision}
\cormark[1]
\cortext[cor1]{Corresponding author.}

\author[1,2]{Shijie Li}[style = Chinese, orcid=0000-0003-0025-760X]
\credit{Writing - Review \& Editing, Supervision}

\author[5]{Yuanchang Liu}[style = Chinese, orcid= 0000-0001-9306-297X]
\credit{Writing - Review \& Editing}
\cormark[1]
\ead{yuanchang.liu@ucl.ac.uk}
\cortext[cor1]{Corresponding author.}

\begin{abstract}
Designing an autonomous precise controller for ships requires accurate and reliable ship models, including the ship dynamic model and actuator model. However, selecting a suitable model for controller design and determining its parameters pose a significant challenge, considering factors such as ship actuation, input constraints, environmental disturbances, and others. This challenge is further amplified for underactuated ships, as obtaining decoupled experiment data is not feasible, and the limited data available may not adequately represent the motion characteristics of ships. To address this issue, we propose a novel model-based parameter estimation approach, called MBPE-LOGO, for underactuated ship motion models. This method combines local optimization and global optimization methods to solve the model identification problem using a dataset generated from real experiments. The effectiveness of the identified model is verified through extensive
comparisons of different trajectories and prediction steps.
\end{abstract}







\maketitle
\section{Introduction} \label{sec1}
\subsection{Background and motivations}

The concept of "Maritime Autonomous Surface Ship" (MASS) has gained attention in recent years as a means of enhancing the intelligence of vessels. Tugs, a type of vessel known for their superior maneuverability compared to traditional large cargo ships, are particularly suitable for implementing intelligent control. Consequently, numerous companies and institutes in different countries have made significant advancements in autonomous tug technology. For instance, the Svitzer Hermod project (2017) \citep{devaraju2018autonomous}, a collaboration between Maersk Line and Rolls Royce, demonstrated the world's first remote-controlled tug in Copenhagen Harbor. The tug was operated remotely by the Remote Operation Center, successfully achieving anchoring and docking with a crew on board. However, this remote operation still required some human presence. Subsequently, in the RECOTUG project (2021) \citep{choi2023review}, Svitzer collaborated with Kongsberg and the American Bureau of Shipping (ABS) to develop a remotely controlled tug capable of performing a full towage operation with all operations controlled from a remote operations center. By remotely controlling tugs, the risk of berthing operations can be reduced and the efficiency of tug services can be enhanced. A similar achievement was made by ASEA Brown Boveri (ABB) and Singaporean shipyard Keppel Offshore $\&$ Marine (Keppel O$\&$M) in Singapore's Port, where they successfully carried out the first remote joystick control of a tug in South Asia. The integration of autonomous technology on the tug enhanced the safety and efficiency of towing operations \citep{choi2023review}. Despite these successes, these projects primarily rely on remote control centers to send instructions for remote control and have not yet achieved complete autonomous control of the vessels themselves.

To achieve autonomous control for tugs, especially in applications such as piloting, escorting, berthing, and unberthing assistance, the establishment of a ship model is essential. Control precision plays a vital role in these scenarios and often requires trajectory tracking control. To design controllers that can perform precise trajectory tracking control, accurate and reliable ship models are required. However, obtaining such ship models can be challenging and often involves the use of Computational Fluid Dynamics (CFD) modeling or sea trials for system identification. These methods require a considerable amount of time and effort, especially when dealing with non-standard vessel types. Furthermore, the established ship maneuvering models, due to their complex parameters, cannot be directly used for ship controller design. Instead, simplified models derived from the ship maneuvering models are commonly used for ship controller design.

The simplified models used for controller design are simplified versions of the ship maneuvering models, particularly the Maneuvering Modeling Group (MMG) model \citep{yasukawa2015introduction,fossen2011handbook}, tailored to different application scenarios. For instance, there are simplified models suitable for low-speed dynamic positioning scenarios and others for normal-speed navigation. The parameters in these models are chosen based on the ship's state, such as speed, with the aim of simplifying the factors that have a relatively minor impact. However, conventional models are often too simplified, such as a model with only six parameters and uncertainty terms, as seen in Remark \ref{rmk:the simplist dynamic model}. Consequently, these simplified models may not fully correspond to the actual vessel, leading to subpar controller performance. To strike a balance between the requirements of controller design and model accuracy, more complex and versatile models have been proposed \citep{saether2019development}.

From the criteria for simplifying ship models, we can deduce that the parameters of a vessel will vary with different ship states, such as speed. This implies that the parameters of the ship model are subject to dynamic changes. However, existing ship parameter model identifications are often based on fixed representative cases, and the represented parameter variations are limited. For example, many studies have used data from Zigzag tests conducted at a fixed speed for ship parameter identification, while ignoring other potential variations in the tests. Therefore, it is worthwhile to explore methods for identifying models that are suitable for varying parameters in different cases. A primary solution for accommodating variable parameters is to employ more complex models and utilize a larger dataset for identification. Complex models are capable of capturing parameter variations, while incorporating a larger dataset into the identification process ensures the validity of more operating conditions. However, this approach brings about new challenges. Not only do complex models make controller design more difficult, but they also hinder effective parameter identification. Moreover, the use of a larger dataset can exacerbate the difficulty of identification to the point where finding the optimal solution may become unfeasible.

In addition, the underactuated ship models face a coupling problem. Each degree of freedom in the ship model exhibits model coupling, resulting in mutual interference of parameters during the identification process. This interference makes it impossible to accurately identify the precise ship model. Motion decoupling is a common solution to mitigate this issue. By performing motion decoupling during the process of designing data acquisition experiments, motion data for individual degrees of freedom can be obtained, allowing for the identification of parameters specific to each degree of freedom. However, if motion decoupling cannot be achieved during data acquisition experiments, accurately identifying motion models becomes a significant challenge.

\subsection{Literature review}

The related works on parameter estimation of autonomous tugs can be divided into two parts: ship dynamic model and experiment design, and parameter estimation methods.

Over the past few decades, ship maneuvering models have been categorized into two types: the Maneuvering Modeling Group (MMG) model \citep{yasukawa2015introduction} and the Abkowitz model \citep{abkowitz1980measurement}. These models primarily aim to examine ship maneuverability. They encompass ship dynamics, propulsion systems, measurement systems, and environmental forces, along with additional features unrelated to control and observer design that affect model accuracy. Consequently, these conventional ship models are too intricate to be directly used in controller design. To address this issue, \cite{fossen2011handbook} has proposed simplified simulation models that facilitate controller design based on different navigation scenarios and tasks. However, simplifying the model relies on prior knowledge of the parameters of the original simulation model, which is not feasible when the maneuvering model is unknown. Therefore, many studies are now focusing on directly identifying parameters using the simplified model. Nevertheless, there are various simplified models available for controller design. Among them, the most commonly applied model for controller design is the most simplified model with only six parameters. Its accuracy is insufficient to accurately reflect the ship's motion trends. Based on the performance comparison between the surge-decoupled model with 22  and the fully-coupled model with 30  in \citep{pedersen2019optimization}, the surge-decoupled model is more likely to be used in real ship controller design compared to other simplified models. The surge-decoupled model, although having a complicated structure to reflect the motion of the actual ship like the fully-coupled model, has fewer parameters, which contributes to the design of the controller. More details and results are presented in Section~\ref{sec5}. After determining the simplified ship dynamic model, the immediate challenge is to determine the parameters. Model-based parameter estimation (MBPE) methods have been proposed, including the least squares (LS) method \citep{huajun2020parameter}, Kalman filtering (KF) \citep{witkowska2020autonomous}, support vector machine (SVM) \citep{xu2012parametric,jiang2021identification}, and optimization methods \citep{du2017simulation,pedersen2019optimization, sirmatel2021modeling}. LS was widely used before the widespread adoption of intelligent algorithms. However, LS is sensitive to outliers in the training sample during the identification process, which can lead to overfitting issues \citep{xu2012parametric}. To address the issues of reliability and serviceability, several improved LS methods have been proposed, such as nonlinear LS, fitting LS, partial LS, and multi-innovation LS \citep{chen2019parameter,sun2013parameter,xie2020parameter}. However, these LS methods are not applicable to the ship controller model due to its nonlinear coupling. KF is also commonly used in ship model identification, and algorithms such as extended KF and untracked KF have emerged \citep{perera2015system,deng2019ukf}. Although these algorithms can provide good performance in ship maneuvering model parameter estimation, they still face issues of parameter drift and dynamic cancellation. SVM, as a supervised learning method, has been utilized in ship model parameter estimation due to its theoretical derivation, generalization capabilities, and global optimality \citep{caccia2000modeling}. However, SVM has a long training time when dealing with large datasets \citep{xue2020system}. Additionally, selecting the kernel function and its parameters is a relatively arbitrary process based on experience \citep{dong2019parameter}. Optimization methods, especially nonlinear optimization, are capable of handling the problem of ship coupling parameters \citep{du2017simulation}. They can also handle large-scale data with noise \citep{sirmatel2021modeling}, which makes them advantageous for ship controller model parameter estimation.

Recent literature shows that optimization methods have been widely used for ship model parameter estimation due to their reliability. In recent years, independent navigation data has enabled the identification of model parameters. For example, in the AutoFerry project, \cite{havdal2017design,pedersen2019optimization} successfully identified the simplified model of a fully actuated ferry ship, named milliAmpere, using the optimization method. By conducting independent navigation experiments for surge, sway, and yaw motion, the parameters of the milliAmpere ship were estimated. Furthermore, control tests using the identified model were demonstrated \citep{bitar2018energy,saether2019development}. However, for underactuated ships, which differ from fully actuated ferries in the AutoFerry project, independent navigation in all three degrees of freedom is not possible. This makes it challenging to obtain motion data for each degree of freedom individually, leading to difficulties in parameter estimation.

\subsection{Original contributions}

Motivated by previous research, this study proposes a comprehensive solution for ship modeling, parameter estimation, and ship motion state validation for a twin thruster tug with underactuated inputs. The surge-decoupled model (22 parameters) of ship motion is chosen based on the milliAmpere project. Other motion control-related models, such as the control allocation model and rotation angle model, have also been identified. To ensure that the identified ship model accurately represents motion characteristics across a wider range of operating conditions, a large dataset from multiple cases is used for the parameter identification process. The proposed solution combines local optimization (LO) and global optimization (GO) to solve the optimization problem effectively. In this study, the main contributions are summarized as follows:

\begin{itemize}
  \item A proposed solution involves the development of a method for modeling underactuated ships and estimating their parameters. The identified models include both the ship's dynamic model and the thruster model. The proposed model aligns with the actual control inputs of the ship, which is advantageous for control purposes.

  \item To solve the parameter estimation problem using a large experimental dataset that considers various ship navigation states, a novel optimization algorithm that combines LO and GO is introduced. Due to the high requirements for initial values in the optimization solution, the estimation of initial values based on empirical formulas helps reduce computational complexity and speed up the process of finding the optimal solution.

  \item The proposed LO method further accelerates the convergence of the algorithm and increases the possibility of finding the optimal solution. The local optimal solution provides the global suboptimal solution and helps the GO find the global optimal solution. The validation results show that the estimated model can have a good fit not only on one step but also on the different time horizons compared to the target trajectories.
  
\end{itemize}

\subsection{Study Organization}

\begin{figure*}[htbp]
  \centering
  \includegraphics[width = 1\textwidth]{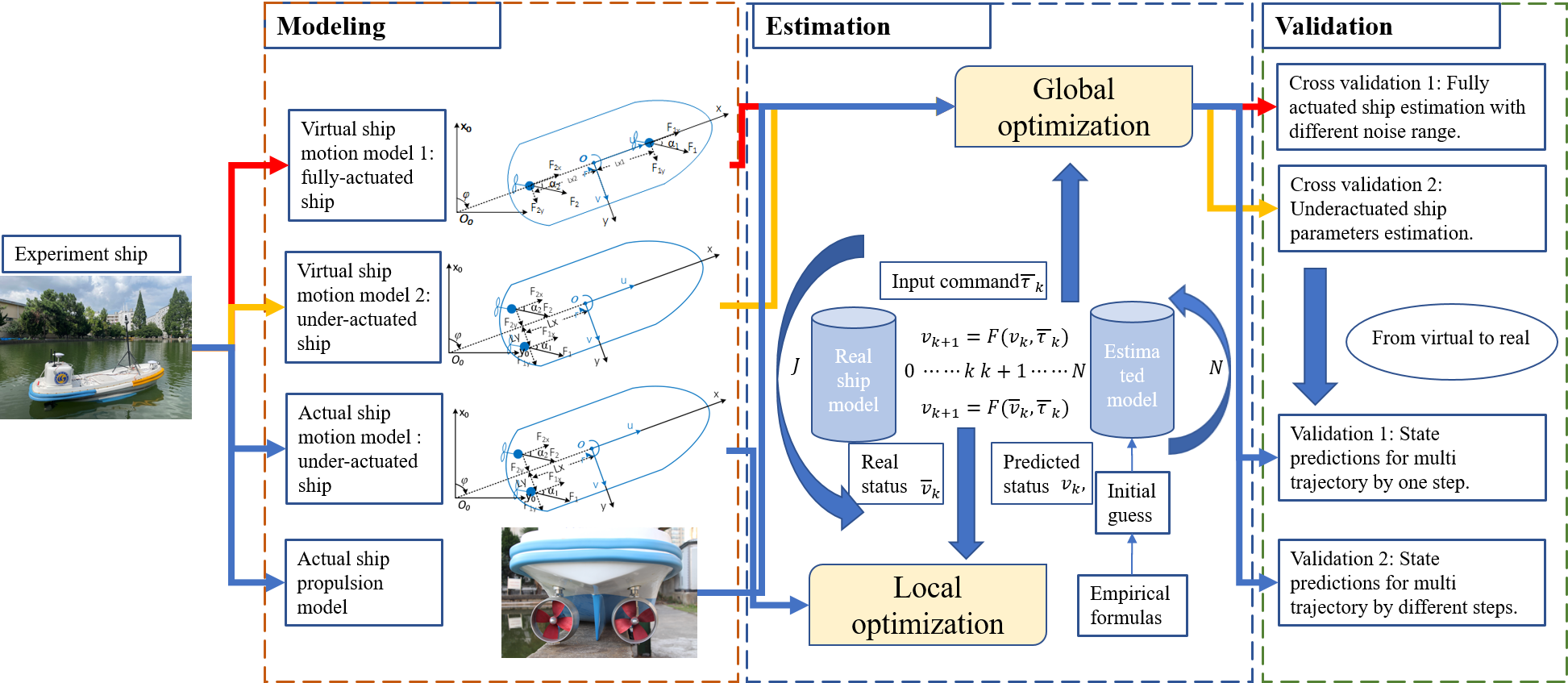}
  \caption{The modeling, estimation, and validation scheme.} 
  \label{f the whole structure of this ship.}
\end{figure*}

The organization of this study is outlined in Figure \ref{f the whole structure of this ship.}, which consists of three main parts: modeling, estimation, and validation. The study uses three main cases for illustration, represented by the red line, yellow line, and blue line. 

The red line represents the verification of the parameter estimation method for the fully actuated ship motion model, which has been proven effective under different disturbance scenarios in a previous study \cite{pedersen2019optimization}. This demonstrates the feasibility of using actual movement data to estimate the real ship parameters and provides reliable parameter estimates under various disturbance scenarios.

Building on this, the yellow line shows the testing of the proposed algorithm on a second virtual ship model with underactuated input. The second ship has the same propulsion system but is configured in an underactuated manner due to the positioning of the propellers. The parameters are assumed to be the same as the first virtual ship. This cross-validation demonstrates the capability of the proposed method to provide an optimal solution for the estimation problem of underactuated ships.

The blue line illustrates the modeling, estimation, and validation procedures of a real tug using the actual ship motion model and propulsion model. The parameters are estimated based on a large amount of estimation data and the complex model structure. To approximate the optimal solution, parameter initialization and the LO method are introduced. Finally, the validations based on real ship movement data confirm the accuracy of the model and the validity of the proposed method.

The structure of this paper is as follows: Section 2 introduces the ship motion model and the actuator model, presenting the modeling and estimation objectives. Section 3 describes the parameter estimation method. Section 4 proposes the validation process for the target ship model estimation, including the introduction of the target ship and the test environment. The whole ship model is determined by combining the thruster model and the ship motion model. Section 5 presents additional real movement data to demonstrate the effectiveness of the proposed whole ship model, with analysis of different prediction steps to evaluate the performance. Finally, Section 6 concludes the paper.

\section{Problem formulation} \label{sec2}

\subsection{Dynamic model of the target ship}

To simplify the surface ship model, it can be described as a 3 degree of freedom (DOF) model, as explained in \cite{fossen2011handbook}:
\begin{figure}[htbp]
  \centering
  \includegraphics[width = 0.45\textwidth]{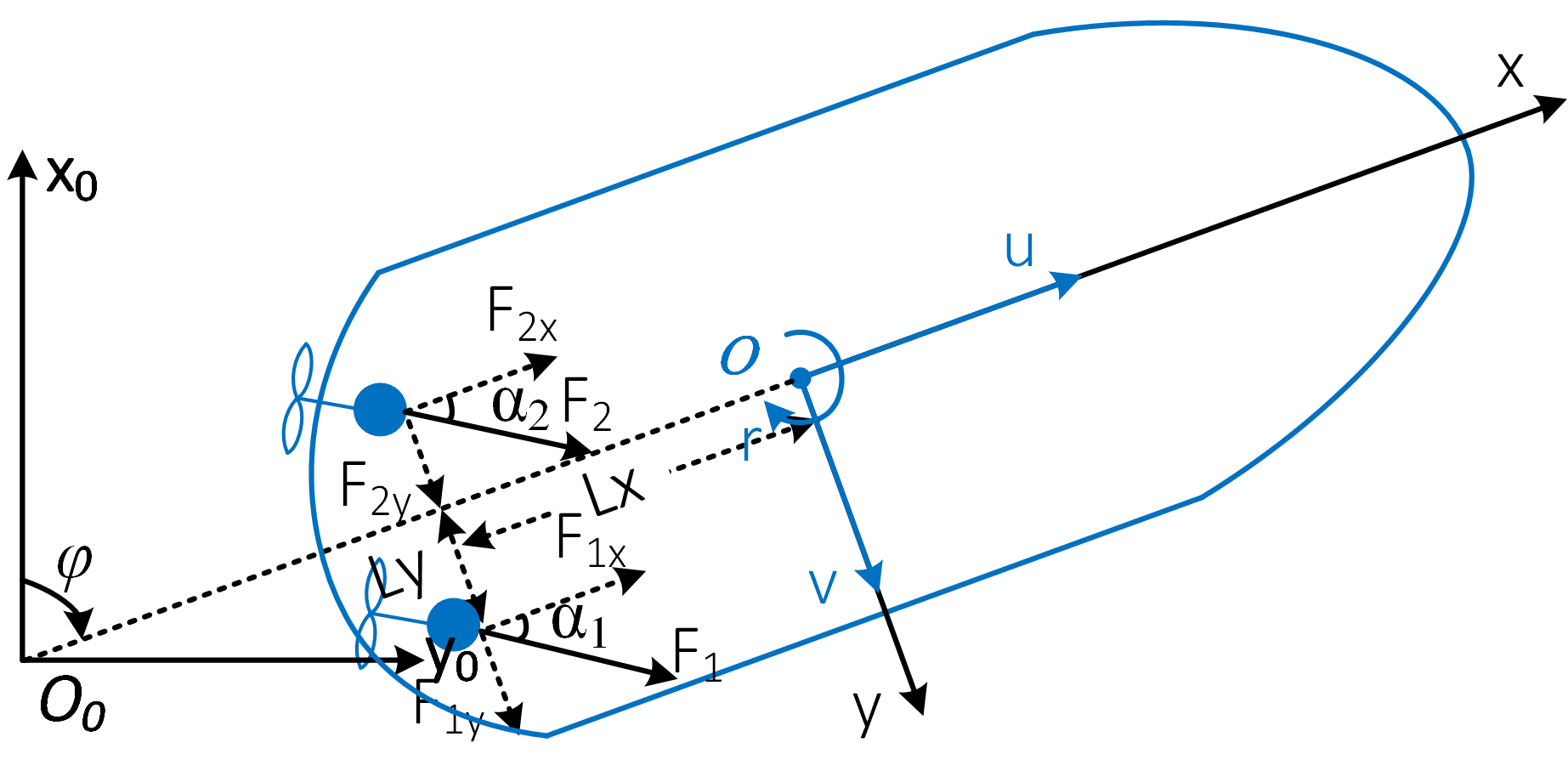}
  \caption{Thruster configuration of the target ship and its coordinate system.}
  \label{f: Thruster configuration of the target ship and its coordinate system.}
\end{figure}

\begin{equation} \label{eq1}
  \begin{cases}
    \boldsymbol{\dot{\eta}} = \boldsymbol{J}(\psi)\boldsymbol{v} \\
    \boldsymbol{M}\boldsymbol{\dot{v}} = -\boldsymbol{C}(\boldsymbol{v})\boldsymbol{v} - \boldsymbol{D(v)}\boldsymbol{v}  + \boldsymbol{\tau},
  \end{cases}
\end{equation}
where the position vector $\boldsymbol{\eta} = [x,y,\psi]^{T}$ represents the ship's position in the earth-fixed coordinate, $x$ and $y$ are the coordinates and $\psi$ is the heading angle. The velocity vector $\boldsymbol{v} = [u,v,r]^{T}$ describes the ship's velocity in the body-fixed coordinate, with $u$ and $v$ denoting the horizontal and vertical components respectively, and $r$ representing the angular speed. The control inputs are denoted by the vector $\boldsymbol{\tau} = [\tau_u, \tau_v ,\tau_r]^T$.

The matrix $\boldsymbol{J}(\psi)$ is a simplified rotation matrix that relates the earth-fixed and body-fixed coordinates. The matrices $\boldsymbol{M}$, $\boldsymbol{C(v)}$, and $\boldsymbol{D(v)}$ correspond to the ship's inertia coefficients, the Coriolis and centripetal matrix, and the damping matrix respectively. To simplify the ship's dynamics, these matrices can be determined based on the surge-decoupled system. The matrices $\boldsymbol{M}$, $\boldsymbol{C(v)}$, and $\boldsymbol{D(v)}$ are defined as follows:

\begin{equation} \label{eq2}
  \begin{aligned}
  \boldsymbol{M} &=\left[\begin{array}{ccc}
  m_{11} & 0 & 0 \\
  0 & m_{22} & m_{23} \\
  0 & m_{32} & m_{33}
  \end{array}\right], \\
  \boldsymbol{C}(\boldsymbol{v}) &=\left[\begin{array}{ccc}
  0 & 0 & c_{13}(v) \\
  0 & 0 &  c_{23}(v)\\
  c_{31}(v) & c_{32}(v) & 0
  \end{array}\right], \\
  \boldsymbol{D}(\boldsymbol{v}) &= \left[\begin{array}{ccc}
  d_{11}(v) & 0 &0 \\
  0 & d_{22}(v) & d_{23}(v) \\
  0 & d_{32}(v) & d_{33}(v)
  \end{array}\right],
\end{aligned}
\end{equation}
where $\boldsymbol{C(v)}$ can be simplified as:
\begin{equation} \label{eq3}
  \begin{aligned}
    c_{13}(\boldsymbol{v}) &= - m_{22} v - m_{23} r \\
    c_{23}(\boldsymbol{v}) &=  m_{11} u \\
    c_{31}(\boldsymbol{v}) &= -c_{13} (v) \\
    c_{32}(\boldsymbol{v}) &= -c_{23} (v),
  \end{aligned}
\end{equation}
and $\boldsymbol{D(v)}$ is determined by:
\begin{equation} \label{eq4}
  \begin{aligned}
    d_{11}(\boldsymbol{v}) &=-X_{u}-X_{|u| u \mid}|u|-X_{u u u} u^{2} \\
    d_{22}(\boldsymbol{v}) &=-Y_{v}-Y_{|v| v}|v|-Y_{|r| v}|r|-Y_{v v v} v^{2} \\
    d_{23}(\boldsymbol{v}) &=-Y_{r}-Y_{|v| r}|v|-Y_{|r| r \mid}|r| \\
    d_{32}(\boldsymbol{v}) &=-N_{v}-N_{|v| v}|v|-N_{|r| v}|r| \\
    d_{33}(\boldsymbol{v}) &=-N_{r}-N_{|v| r}|v|-N_{|r| r}|r|-N_{r r r} r^{2},
  \end{aligned}
\end{equation}
where the ship parameters are therefore determined:
\begin{equation} \label{eq5}
  \begin{split}
      \boldsymbol{P} =[m_{11}, m_{22}, m_{23}, m_{32}, m_{33}, X_{u}, X_{|u| u}, X_{u u u}, Y_{v}, Y_{|v| v}, Y_{v v v}, \\
      Y_{|r| v}, Y_{r}, Y_{|v| r}, Y_{|r| r}, N_{v}, N_{|v| v}, N_{|r| v}, N_{r}, N_{|r| r}, N_{r r r}, N_{|v| r}].
  \end{split}
\end{equation}

The ship model can be selected using the fully coupled system described in \citep{pedersen2019optimization}, which showed that the differences between the surge-decoupled system and the fully-coupled system are negligible. Additionally, the surge-decoupled system has fewer unknown parameters for estimation, making it a suitable choice for the target ship model in this research. Furthermore, compared to the simplest model in Remark \ref{eq simplified dynamic model}, the surge-decoupled model is more complex with a greater number of parameters, allowing it to capture more dynamic characteristics. From Equations \ref{eq1} and \ref{eq2}, it can be observed that the input torques have values in all three degrees of freedom (DOFs) and are coupled to each other, resulting in a matrix that is not full-rank. Consequently, the ship dynamic model is under-actuated.

\begin{remark} \label{rmk:the simplist dynamic model}
With a similar structure to Equation \eqref{eq1}, the most simplified dynamic model only includes six parameters $(m_{11},m_{22}, m_{33},d_{11},d_{22},d_{33})$ in Equation \ref{eq simplified dynamic model}, as mentioned in \citep{chen2020leader,yuan2021leader}. The parameter count of this model is lower than the selected model due to the exclusion of additional factors, which may result in reduced accuracy. In the following section, this model will be used for comparison purposes.
  \begin{equation} \label{eq simplified dynamic model}
    \begin{aligned}
    \boldsymbol{M} &=\left[\begin{array}{ccc}
    m_{11} & 0 & 0 \\
    0 & m_{22} & 0 \\
    0 & 0 & m_{33}
    \end{array}\right], \\
    \boldsymbol{C}(\boldsymbol{v}) &=\left[\begin{array}{ccc}
    0 & 0 & -m_{22}v \\
    0 & 0 &  m_{11}u\\
    m_{22}v & -m_{11}u & 0
    \end{array}\right], \\
    \boldsymbol{D}(\boldsymbol{v}) &= \left[\begin{array}{ccc}
    d_{11} & 0 &0 \\
    0 & d_{22} & 0 \\
    0 & 0 & d_{33}
    \end{array}\right],
  \end{aligned}
  \end{equation}
\end{remark}

\subsection{The thruster model}
\subsubsection{Thrust transformation}
Typically, the inputs of the ship dynamic model are torques rather than actual control commands, such as rudder angles and shaft speed commands. The ship's propulsion system configuration determines whether the ship is fully actuated or underactuated. The target ship's thruster configuration is illustrated in \ref{f: Thruster configuration of the target ship and its coordinate system.}. According to \citep{fossen2011handbook}, the azimuth thruster and its location can be described as:
\begin{equation} \label{eq6}
  \boldsymbol{\tau} = \boldsymbol{T}(\boldsymbol{\alpha})\boldsymbol{F},
\end{equation}
where $\boldsymbol{F} = [f_1,f_2]^{\rm{T}}$ represents the forces generated by the thruster speed, and the thruster configuration matrix is calculated as:
\begin{equation} \label{eq7}
  \begin{aligned}
  T(\alpha) &= \left[\begin{array}{cc}
  \cos(\alpha_1) & \cos(\alpha_2)  \\
  \sin(\alpha_1) & \sin(\alpha_2)  \\
  l_{x1} \sin(\alpha_1) + l_{y1} \cos(\alpha_1) & l_{x2} \sin(\alpha_2) + l_{y2} \cos(\alpha_2)
  \end{array} \right].
  \end{aligned}
\end{equation}

As depicted in Figure \ref{f: Thruster configuration of the target ship and its coordinate system.}, the azimuth thruster configuration of the target ship and its coordinate system exhibit a clockwise rotation along the positive BODY x-axis in accordance with the BODY fixed coordinate direction. Consequently, the thruster force can be divided into the x-axis and y-axis directions. For the torque in the yaw direction, the determination of $l_x$ and $l_y$ (positive or negative) depends on the position from $O$. Both azimuth rotation angle units are represented in radians. Therefore, the direction of the results corresponds to Equation \ref{eq1}.

\subsubsection{The thruster force model}
The force exerted on the shaft by the propeller is highly nonlinear and influenced by various factors such as propeller shape, relative velocity of water flow, and pressure fluctuations within the hull's wake \citep{overaas2020dynamic}. To simplify this dynamic relationship, a polynomial fitting approach has been utilized.

\subsubsection{The azimuth thruster rotation angle model}
\begin{figure*}[htbp]
  \centering
  \includegraphics[width = 1\textwidth]{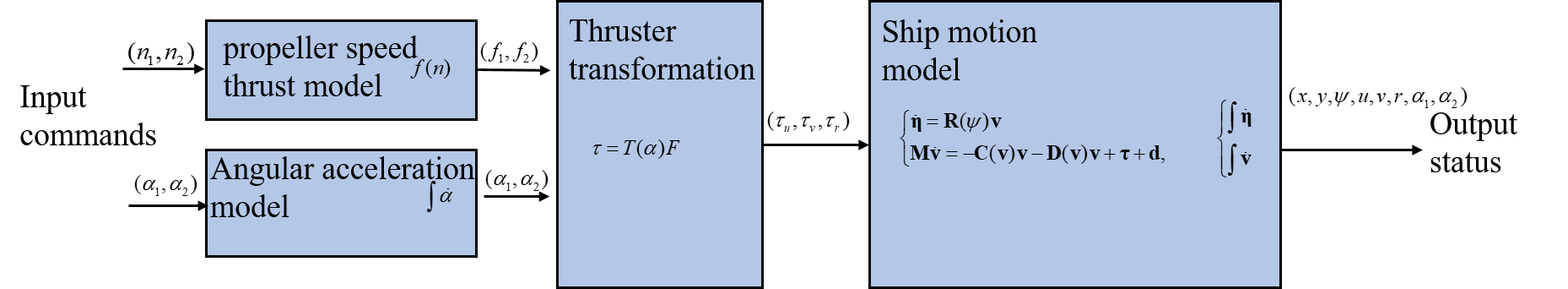}
  \caption{The whole structure of the ship motion model.}
  \label{f the whole ship motion model}
\end{figure*}

Based on the performance of the thruster angle rotation, the dynamics of the azimuth angle demonstrate an S-shaped behavior. Hence, the model can be described as follows:
\begin{equation}
  \dot{\alpha} = K_{\alpha} \frac{(\alpha_d - \alpha)}{\sqrt{(\alpha_d - \alpha)^2 + \epsilon^2 }},
\end{equation}
where $\alpha$ represents the current azimuth angle, $\alpha_d$ is the command or desired azimuth angle, and $K_{\alpha}$ and $\epsilon $ are constant parameters that describe the rotation change speed. Consequently, the entire motion model of the ship, encompassing both the thruster and the ship dynamic models, can be determined (see Figure \ref{f the whole ship motion model}).

\section{The parameter estimation method of the ship motion model} \label{sec3}
In this section, we introduce the methodology for initializing ship parameters, along with the LO method and GO method for model-based parameter estimation. Additionally, two cross-validations are proposed for verification purposes. Cross-validation 1 seeks to confirm the effectiveness of the GO methods in estimating fully actuated ship parameters under various perturbation conditions, while cross-validation 2 further verifies their effectiveness in estimating under-actuated ship parameters.

\subsection{Ship parameters initialization}
Careful consideration must be given to the initial values for parameters, particularly in scenarios that involve complex and diverse parameters and coupled models. For the surge-decoupled ship model, certain parameters can be estimated using empirical formulas based on their definitions. In this study, the initialization of the ship parameters is determined using empirical formulas according to the ship status and configurations. However, only some parameters can be easily calculated, including $[m_{11},m_{22},m_{23},m_{32},m_{33}]$. The other parameters are set with initial values of zero.

$\boldsymbol{M} $ is defined as follows:
\begin{equation}
  \boldsymbol{M}:=\left[\begin{array}{ccc}
  m-X_{\dot{u}} & 0 & 0 \\
  0 & m-Y_{\dot{v}} & m x_g-Y_{\dot{r}} \\
  0 & m x_g-N_{\dot{v}} & I_z-N_{\dot{r}}
  \end{array}\right],
\end{equation}
where $ m $ represents the total mass, $ X_{\dot{u}}, Y_{\dot{v}}, Y_{\dot{r}}, N_{\dot{v}}, N_{\dot{r}} $ are the hydrodynamic derivatives, $ x_g $ is the distance from the gravity center to the original center. 

Referred to in Section 2.4.1 in \cite{fossen1999guidance}, employing strip theory allows for the estimation of some hydrodynamic derivatives as follows:
\begin{equation} \label{eq10}
  \begin{cases}
  \begin{aligned}
    X_{\dot{u}} &= -0.05 m \\
    Y_{\dot{v}} &= -0.5\rho D^2 L \\
    N_{\dot{r}} &= \frac{1}{24}(0.1 m B^2 + \rho \pi D^2 L^3) \\
    I_y &= I_z = \frac{4}{15} \pi \rho a b^2 (a^2 + b^2) ,
  \end{aligned}
  \end{cases}
\end{equation}
where $L$ is the ship length, $B$ is the ship width, $D$ is the hull draft, $\rho$ is the water density, $a,b$ are the half of the $L,B$ separately.  Usually, $ x_g $  is supposed to be 0,  $Y_{\dot{r}}, N_{\dot{v}}$ are very small. Thus, the corresponding parameters are determined, and the specific values of this experiment ship are demonstrated in Equation \ref{eq15}. 

\begin{remark}
  By using empirical formulas, it is possible to obtain estimates of the corresponding parameters. However, these estimates are not true values, and they do not remain completely consistent in the estimation process.
\end{remark}

\subsection{Optimization algorithm for ship parameters estimation}

The ship parameters estimation problem can be formulated based on Equation \ref{eq1}. After determining the model structures and initial values, the estimation of parameters is performed using the optimization method. However, since the optimization control problem can easily get trapped in local optima, it is necessary to introduce multiple constraints to ensure solution conditions. To accelerate the search for the global optimal solution, the LO structure is utilized to find a suboptimal solution. This allows the GO method to provide the global optimal solution faster and more efficiently. The optimization solver chosen for this study is IPOPT on CasADi \cite{Andersson2019}, which is effective for solving large-scale nonlinear optimization problems in continuous systems.

\subsubsection{The GO method}
Under the assumption that the transformation between NED coordinates and BODY coordinates can be simplified as $\boldsymbol{J}({\boldsymbol{\psi}})$, the identification problem can be reduced to finding the relationship involving $\dot{\boldsymbol{v}}$. By collecting the ship input data $\boldsymbol{\tau}$ and output data $\boldsymbol{v}$, the problem formulation for the GO method is as follows:
\begin{equation}
  \begin{array}{ll}
  \min\limits_{{\boldsymbol{w}}} & \phi(\boldsymbol{w}, \overline{\boldsymbol{w}}) \\
  \text { s.t } & \boldsymbol{g}(\boldsymbol{w})=\mathbf{0} \\
  & h(\boldsymbol{w}) \leq 0 \\
  & \boldsymbol{w}_{l b} \leq \boldsymbol{w} \leq \boldsymbol{w}_{u b} \\
  & \boldsymbol{w}(0)=\overline{\boldsymbol{w}}(0),
  \end{array}
\end{equation}
where
\begin{equation}
  \begin{aligned}
  \phi(\boldsymbol{w}, \overline{\boldsymbol{w}}) &=\sum_{j=1}^J \sum_{n=1}^{N_j}\left(\boldsymbol{v}_{j n}-\overline{\boldsymbol{v}}_{j n}\right)^T \boldsymbol{W}_j\left(\boldsymbol{v}_{j n}-\overline{\boldsymbol{v}}_{j n}\right)+\lambda \boldsymbol{R}(\boldsymbol{P}) \\
  \boldsymbol{w} &=\left[\boldsymbol{v}_{10}^T, \boldsymbol{v}_{11}^T, \ldots, \boldsymbol{v}_{1 N}^T, \boldsymbol{v}_{20}^T, \boldsymbol{v}_{21}^T, \ldots, \boldsymbol{v}_{J N_J}^T, P^T\right]^T \\
  \overline{\boldsymbol{w}} &=\left[\overline{\boldsymbol{v}}_{10}^T, \overline{\boldsymbol{v}}_{11}^T, \ldots \overline{\boldsymbol{v}}_{1 N_1}^T, \overline{\boldsymbol{v}}_{20}^T, \overline{\boldsymbol{v}}_{21}^T, \ldots \overline{\boldsymbol{v}}_{J N_J}^T\right]^T,
  \end{aligned}
\end{equation}
where $J$ denotes the number of maneuvers, and $N_j$ represents the number of samples of each maneuver. $W$ stands for the weight of the maneuvers. $\boldsymbol{w}$ corresponds to the decision vector which is subject to change during the estimation process. $\overline{\boldsymbol{w}}$ represents the experimental data. $W$ varies depending on the specific maneuver. $\lambda \boldsymbol{R}(\boldsymbol{P})$ restricts the size of parameters through a ridge regression approach. $\lambda$ is a constant that determines the tradeoff between bias and variance. $\boldsymbol{R}(\boldsymbol{P})$ denotes the norm $\|\boldsymbol{P}  \|_2  $. The equality constraints are given by:
\begin{equation} \label{eq13}
  \begin{aligned}
    \boldsymbol{v}_{11} &= \boldsymbol{F}(\boldsymbol{v}_{10},\overline{\boldsymbol{\tau}}_{10}) \\
    g(\boldsymbol{w}) &=\left[\begin{array}{c}
      \boldsymbol{F}\left(\boldsymbol{v}_{10}, \overline{\boldsymbol{\tau}}_{10}\right)-\boldsymbol{v}_{11} \\
      \ldots \\
      \boldsymbol{F}\left(\boldsymbol{v}_{1 N_1-1}, \overline{\boldsymbol{\tau}}_{1 N_1-1}\right)-\boldsymbol{v}_{1 N_1} \\
      \boldsymbol{F}\left(\boldsymbol{v}_{20}, \overline{\boldsymbol{\tau}}_{20}\right)-\boldsymbol{v}_{21} \\
      \ldots \\
      \boldsymbol{F}\left(\boldsymbol{v}_{J N_J-1}, \overline{\boldsymbol{\tau}}_{J N_J-1}\right)-\boldsymbol{v}_{J N_J}
      \end{array}\right]=\mathbf{0} ,\\
  \end{aligned}
\end{equation}
where $\boldsymbol{F}$ denotes the calculation function for ship speed according to Equation \ref{eq1} using 4th order Runge-Kutta (RK4). In this case, the GO method is employed and uses the states computed from the estimated model at each step of the optimization process. It then solves the error minimization constraint between the calculated states and the collected true states.

\subsubsection{The combination of the GO and LO method}
To modify the Model-Based Parameter Estimation (MBPE) process, the LO method considering the convex approximation (CA) is introduced:
\begin{equation} \label{eq14}
  \begin{aligned}
    g(\boldsymbol{w}) &=\left[\begin{array}{c}
      \boldsymbol{F}\left(\overline{\boldsymbol{v}}_{10}, \overline{\boldsymbol{\tau}}_{10}\right)-\boldsymbol{v}_{11} \\
      \ldots \\
      \boldsymbol{F}\left(\overline{\boldsymbol{v}}_{1 N_1-1}, \overline{\boldsymbol{\tau}}_{1 N_1-1}\right)-\boldsymbol{v}_{1 N_1} \\
      \boldsymbol{F}\left(\overline{\boldsymbol{v}}_{20}, \overline{\boldsymbol{\tau}}_{20}\right)-\boldsymbol{v}_{21} \\
      \ldots \\
      \boldsymbol{F}\left(\overline{\boldsymbol{v}}_{J N_J-1}, \overline{\boldsymbol{\tau}}_{J N_J-1}\right)-\boldsymbol{v}_{J N_J}
      \end{array}\right]=\mathbf{0} ,\\
  \end{aligned}
\end{equation}
where $\overline{\boldsymbol{v}}$ denotes the experimental data. The main difference between the GO and LO methods lies in the selection of calculation inputs. The LO method selects experimental data rather than calculation data to quickly determine the estimated parameter values. However, the other steps involved in parameter estimation remain the same. Detailed algorithmic information can be found in Figure~\ref{f the whole structure of this ship.}.

If the predicted model is 100\% consistent with the true model, the estimated model should have the same state values as the true model when given the same initial values and continuous control inputs. However, achieving this idealized result of model identification and the optimal solution of the GO method is difficult due to external disturbances and modeling errors, among other factors. The difference between LO and GO lies in whether the actual ship position information can be continuously obtained during every optimization calculation process. The LO method can acquire this information, resulting in a faster calculation speed and lower computing power consumption, and it is guaranteed to have a solution. On the other hand, the GO method may not be able to find a solution due to the heavier calculation effort. Therefore, one of the main innovations in this study is to choose the LO result as the reference result for the GO method. For additional detailed analysis, please refer to \cite{bonilla2010convexity, Andersson2019, sirmatel2021modeling}.

Due to the large amount of data, importing all the experimental data into the GO process may lead to no solution. To address this issue, a method of combining the LO and GO methods is proposed, gradually increasing the quantity of input data. The algorithm for this method is presented in Algorithm~\ref{alg1}.

\begin{algorithm}[h]
	\renewcommand{\algorithmicrequire}{\textbf{Input:}}
	\renewcommand{\algorithmicensure}{\textbf{Output:}}
	\caption{The combination with LO and GO method.}
	\label{alg1}
	\begin{algorithmic}[1]
        \REQUIRE The model structure, input data, LO method, and GO method.
		\STATE Initialization: $ n $
        \WHILE{$n != n_{max}$}
        \STATE Select the input data according to $C_n^{n_{max}}$ maneuver.
        \STATE Using the GO method to get the result $P_{GO}$.
            \IF{No result}
            \STATE Using the LO method to get the result $P_{LO}$.
            \STATE Using the $P_{LO}$ as the initial guess in GO process to get $P_{GO}$.
            \ENDIF
        \STATE $n = n+1 $
        \ENDWHILE
		\ENSURE The result $P_{GO}$.
	\end{algorithmic}  
\end{algorithm}
where $n$ represents the number of maneuver cases, $n_{max}$ represents the total number of maneuver cases. The variables $P_{GO}$ and $P_{LO}$ represent the parameter estimation results calculated by the GO and LO methods, respectively.  In cases where the GO process is unable to compute the optimal solution, the results obtained from the LO method can be used as initial values instead. This ensures that even with a large dataset, an optimal solution can still be achieved.

\section{The target ship model identification}
In this section, a real ship model is introduced to verify the proposed parameter estimation method. The surge decoupled model is used during the LO and GO processes. Moreover, the relationship between the thruster speeds and the forces has been identified via polynomial regression. Therefore, the whole structure of the ship motion model from the actual inputs and actual outputs can be determined.

\subsection{The research object: Qiuxin No.5 tug }
\begin{figure}[htbp]
  \centering
  \includegraphics[width = 0.45\textwidth]{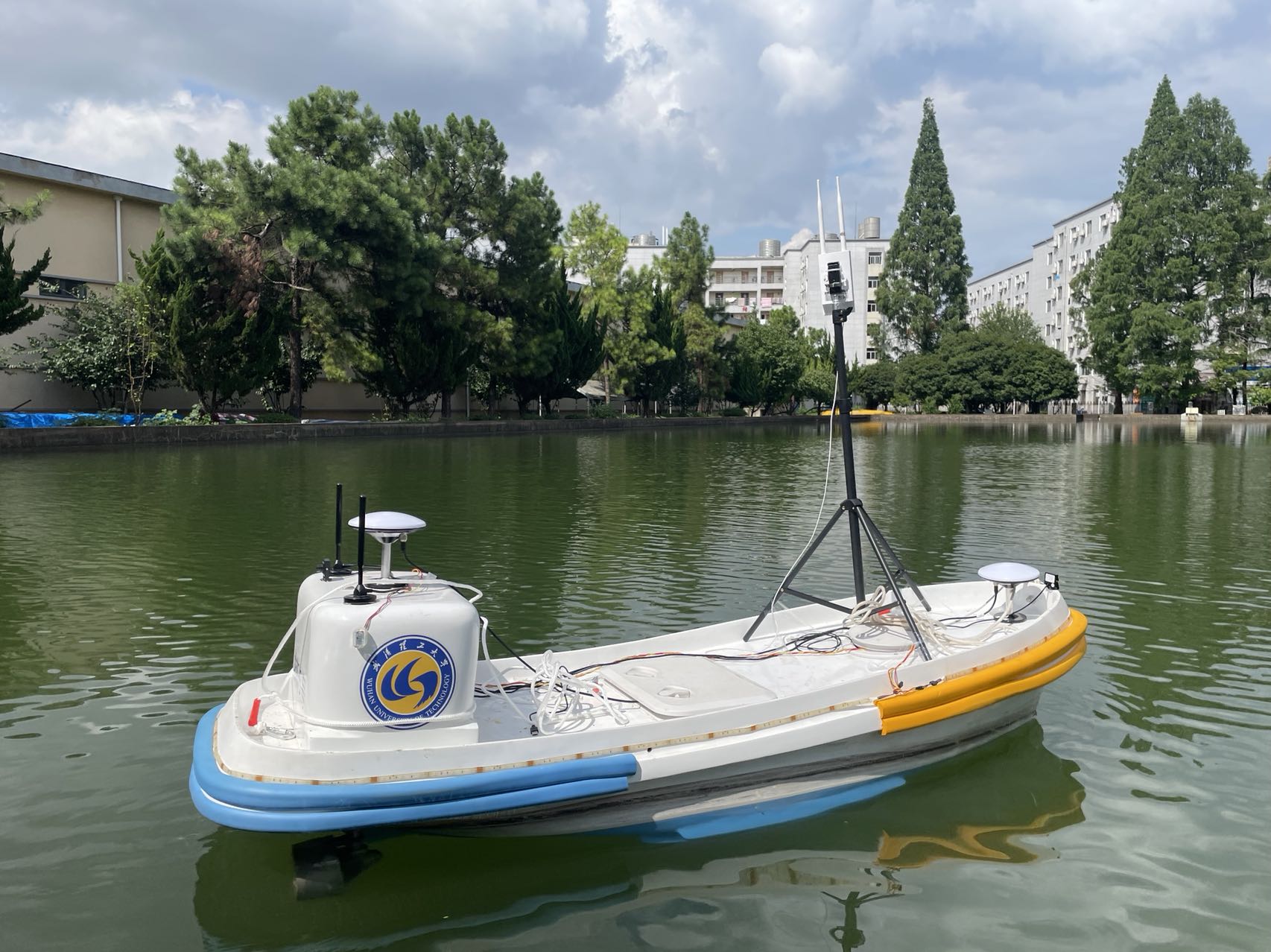}
  \caption{The navigation of the Qiuxin No.5.}
  \label{f:Qiuxin No.5}
\end{figure}

\begin{figure}[htbp]
  \centering
  \includegraphics[width = 0.45\textwidth]{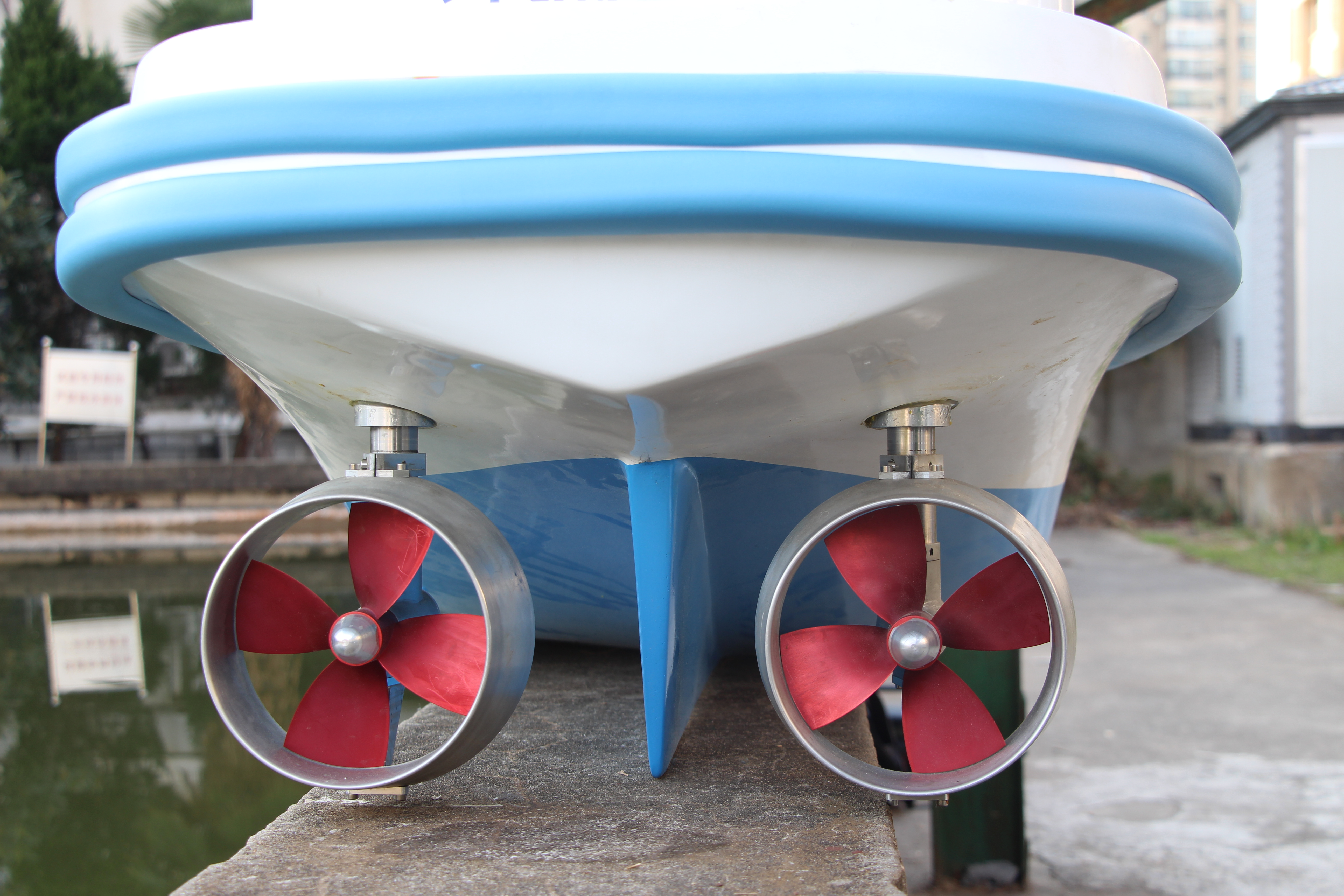}
  \caption{The azimuth thruster configuration of Qiuxin No.5.}
  \label{f:azimuth thrusters}
\end{figure}

The Qiuxin No.5 tug (shown in Figure~\ref{f:Qiuxin No.5}) is an autonomous tug ship designed for autonomous testing. It is a 1:20 scaled tug commonly used in some Chinese ports. The tug is symmetrical with two azimuth thrusters, as shown in Figure \ref{f:azimuth thrusters}. The azimuth thruster is located aft of the vessel, 0.8m from the CO and 0.163m from the centerline, as shown in Figure \ref{f: Thruster configuration of the target ship and its coordinate system.}. The vessel is fully electric and has a battery bank of 48v, 28Ah. The specifications for the Qiuxin No.5 tug are present in Table \ref{t:Qiuxin No.5 parameters}. This tug has a Global Positioning System (GPS) sensor with Real-time kinematic positioning (RTK), and an inertial measurement unit (IMU) sensor to obtain accurate position and speed information. The GPS sensor provides GPS coordinates, East, North, and Up (ENU) positions, and speed information. The IMU sensor provides 3-axis acceleration in BODY fixed coordinates and 3-axis angular velocity around 3-axis acceleration. Finally, all data frequencies are unified to 0.2 Hz.
\begin{table}[htbp]
  \centering
  \caption{Qiuxin No.5 parameters.} \label{t:Qiuxin No.5 parameters}
    \begin{tabular}{ccc}
      \hline
      Parameters & Values  & Unit \\
      \hline
      $dispv$ & 0.1876  & $m^3$ \\
      $\rho$ & 1000 & $kg/m^3$ \\
      $m$ & 187.600 & $kg$ \\
      $D$ & 0.2485    & $m$\\
      $B$ & 0.6952    & $m$ \\
      $L$ & 2.152    & $m$\\
      \hline
    \end{tabular}
\end{table}

\begin{figure}[htbp]
  \centering
  \includegraphics[width = 0.45\textwidth]{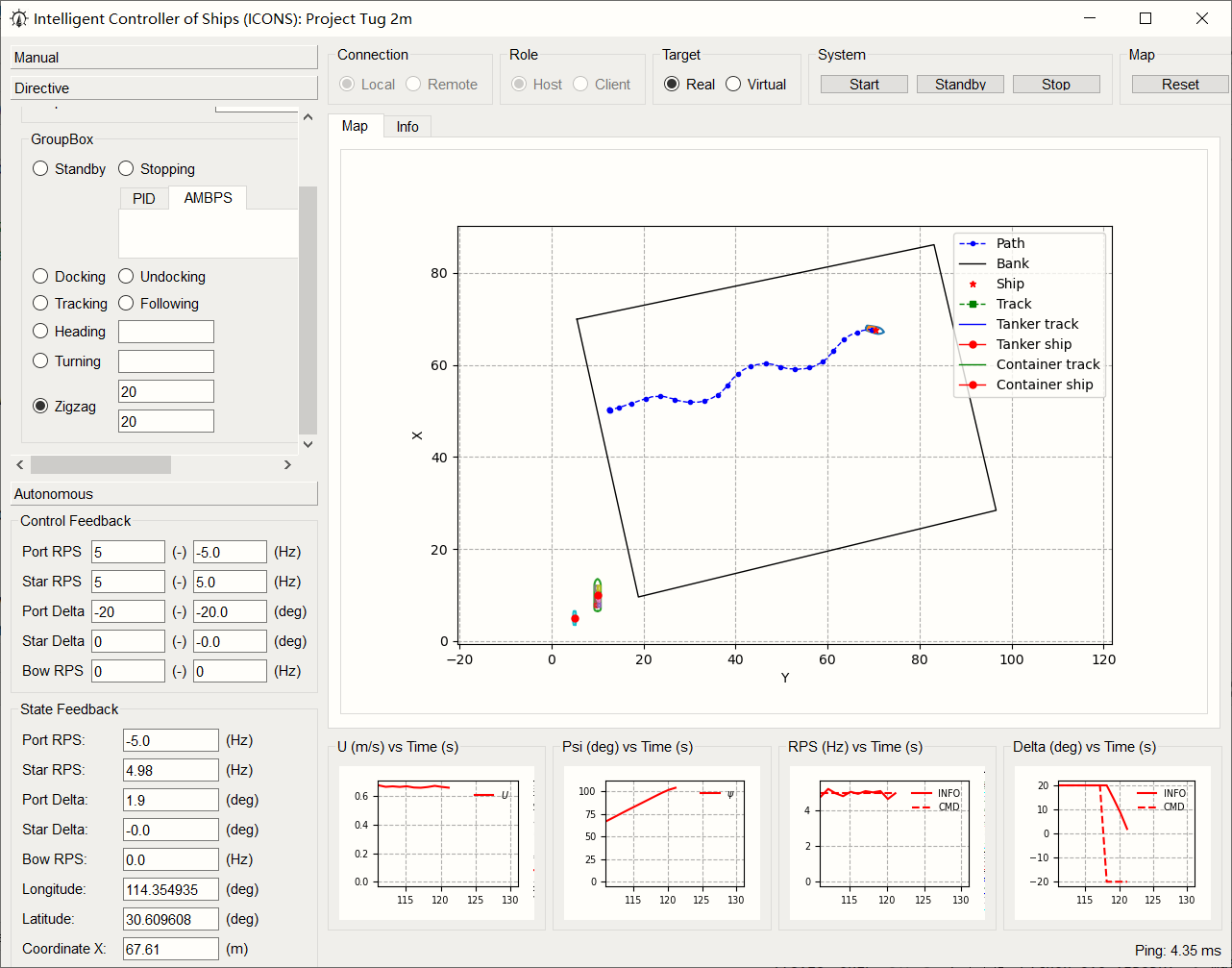}
  \caption{+20+20 zigzag experiment.}
  \label{f:navigation}
\end{figure}

As mentioned in Section \ref{sec3}, for under-actuated ships, the experimental data of straight lines, turning circles, and zigzags can be used for ship model parameter estimation. Thus, in the 60*80m pool, these tests are carried out under low disturbances as shown in Figure \ref{f:navigation} by using software referred to in \cite{you2020integrated} and based on the structure outlined in \cite{yan2019applying}.

\subsection{The surge-decoupled model identification}
\begin{figure*}[htbp]
  \centering
  \includegraphics[width = 1\textwidth]{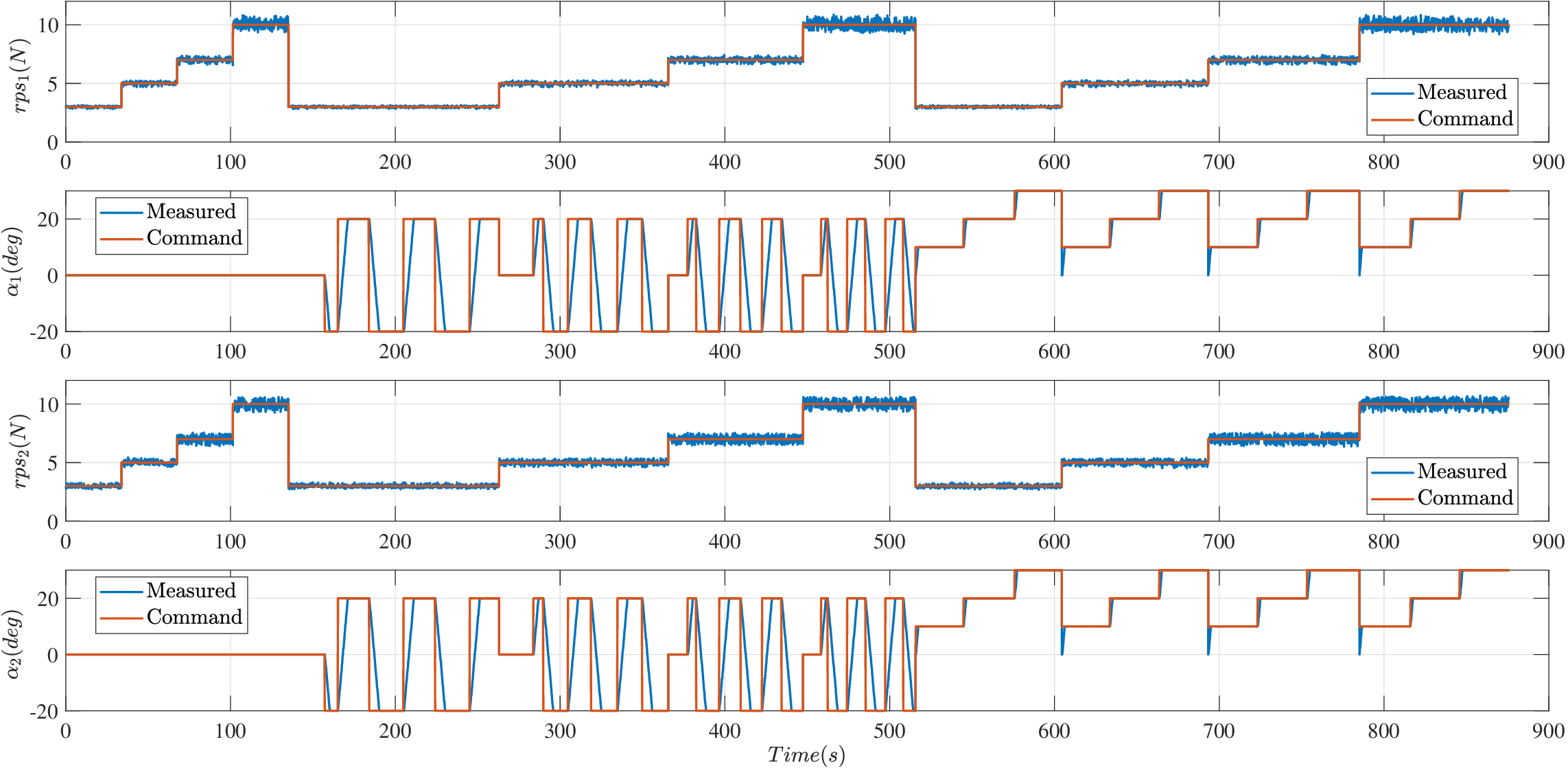}
  \caption{The input command feedback data of the designed 12 maneuvers.}
  \label{f real feedback command}
\end{figure*}

The surge-decoupled model is identified using a set of 12 maneuvers, as shown in Table 2. These maneuvers include 4 straight lines with constant thruster shaft speeds of 3, 5, 7, and 10 revolutions per second (RPS), 4 +20+20 zigzag tests with the corresponding thruster speeds, and 4 turning circles with the corresponding thruster speeds. In order to increase maneuverability, the turning circle is varied by changing the azimuth angles from 10 deg to 30 deg. The real input feedback data collected from the sensors in the ship are shown in Figure \ref{f real feedback command}, where the low-speed data is smooth and the high-speed data fluctuates slightly. Consequently, the input torques also exhibit some fluctuations, as shown in Figure \ref{f real feedback command}. It is worth noting that the high shaft speed fluctuations can make it challenging to obtain an optimal solution for parameter estimation.

\begin{table*}[h]
    \caption{The designed 12 maneuvers for the under-actuated target ship parameter estimation.}
    \label{tab:12 maneuvers}
    \centering
    \begin{tabular}{cc}
        \hline
        12 maneuvers & Cases\\
        \hline
        1st to 4th & Straight line (u-direction) with 3,5,7,10 RPS shaft speed\\
        5th to 8th & +20+20 Zigzag movements with 3,5,7,10 RPS shaft speed \\
        8th to 12th & Turning circle movements with 10, 20, 30 deg rotation angle  \\
        \hline
    \end{tabular}
\end{table*}

\begin{figure*}[htbp]
  \centering
  \includegraphics[width = 1\textwidth]{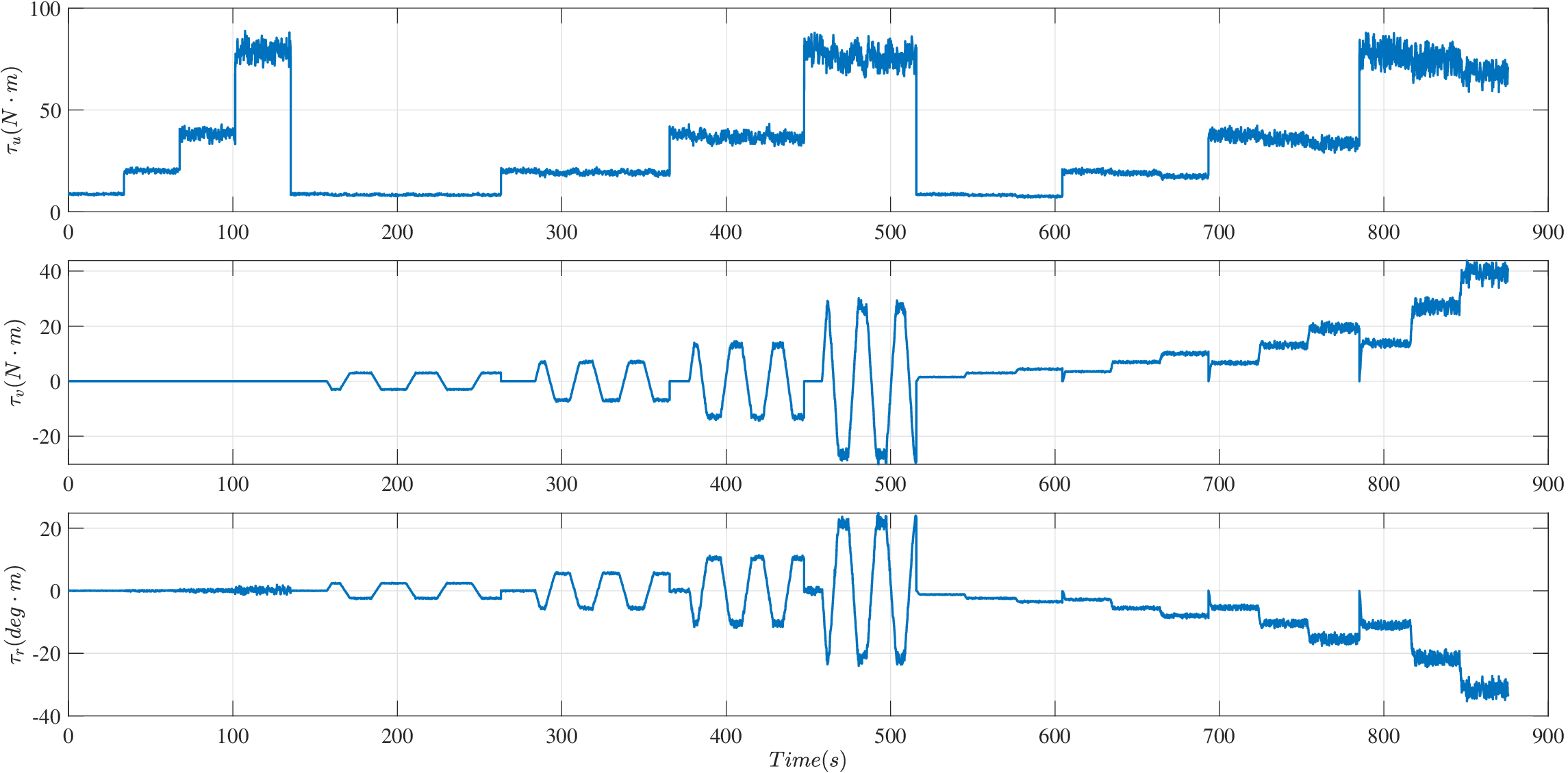}
  \caption{The calculated input torques from the real input in Figure \ref{f real feedback command}.}
  \label{f calculated input torques}
\end{figure*}

\begin{remark}
  In zigzag experiments, different symbols and values represent different settings. For example, in the +20+20 zigzag experiment, the former +20 is the maximum difference between the actual heading angle and the initial heading angle during the motion, and in the second +20, the positive sign represents the initial rotation angle direction, with the value 20 representing the maximum rotation angle degree.
\end{remark}

\begin{remark}
The literature \cite{jiang2022identification, carrica2013turn} only uses the +20+20 zigzag data to determine the ship motion model. However, this approach is not accurate. The ship motion model is very complex and may contain time-varying parameters and more parameters than the simplified model. When the model structure is fixed, using more cases helps to reflect a wider range of motion characteristics. The results show that the identified model is suitable for all cases.
\end{remark}

Considering the ship modeling and ignoring the influence of the delay of the controller, the feedback values (also measured values and supposed as the true values) are taken as the input data to calculate the input torques. As for the output data, due to the lack of speed sensors, $[u,v,r]$ data is not able to be obtained directly. However, it can be converted from data obtained from IMU sensors. The conversion formulas are selected as $\boldsymbol{v} = \boldsymbol{J}(\psi)\dot{\eta}$, where $\dot{\eta}$ comes from the NED position from GPS data.

According to the analysis in Section \ref{sec3}, the estimation method is greatly affected by the initial values, which can be calculated along with the empirical formulas in Equation \ref{eq10} and ship configurations in Table \ref{t:Qiuxin No.5 parameters}.
\begin{equation} \label{eq15}
  \begin{cases}
  \begin{aligned}
    X_{\dot{u}} &= -0.05 m = -9.38 \\
    Y_{\dot{v}} &= -0.5\rho D^2 L = -66.4454\\
    N_{\dot{r}} &= -\frac{1}{24}(0.1 m B^2 + \rho \pi D^2 L^3) =-80.9376\\
    I_y &= I_z = \frac{4}{15} \pi \rho a b^2 (a^2 + b^2)=139.2598 .
  \end{aligned}
  \end{cases}
\end{equation}

Hence, the parameter initialization for $\boldsymbol{M}$ is introduced as:
\begin{equation}
  \boldsymbol{M} =\left[\begin{array}{ccc}
  196.98 & 0 & 0 \\
  0 & 254.0454 & 0 \\
  0 & 0 & 220.1974
  \end{array}\right].
\end{equation}

As the initialization for $\boldsymbol{M}$ has been determined, we make the following assumptions for parameter estimation. We select the inequality of the inertia mass matrix as:
\begin{equation}
  \begin{aligned}
    m_{11} >0; m_{22} >0;  m_{33} >0;
  \end{aligned}
\end{equation}

Also, the optimization algorithm may have no solution due to the complexity of the large-scale data, as the total data sample is 4390. To overcome this and simplify the estimation process, the inequalities are introduced according to Theorem 5.2 in \cite{fossen1999guidance,pedersen2019optimization}.
\begin{equation}
  h(\boldsymbol{w}) = \left[
    \begin{array}[]{c}
      X_u \\
      X_u(Y_v N_r - Y_r N_v)
    \end{array}
  \right] > 0
\end{equation}
which can limit the range of the optimization solution, and speed up the optimization process. These inequalities are under the assumption that the ship is dynamically stable in straight-line motion, where the target ship meets this condition.

\begin{figure*}[htbp]
  \centering
  \includegraphics[width = 1\textwidth]{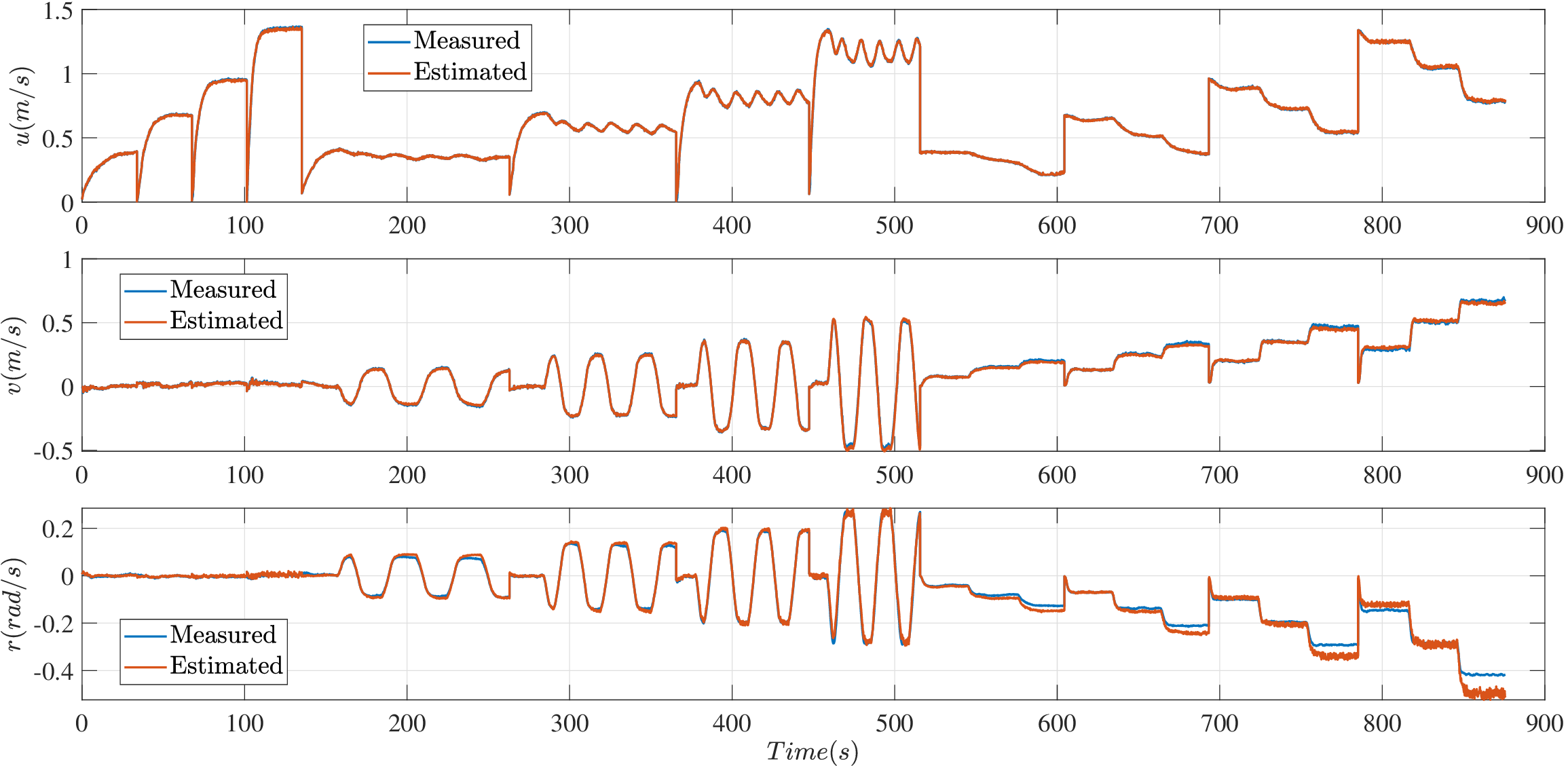}
  \caption{The data validation from the estimated result $P_{LO}$ by using $C_{12}^{12}$ maneuver data.}
  \label{f CA validation}
\end{figure*}

\begin{figure*}[htbp]
  \centering
  \includegraphics[width = 1\textwidth]{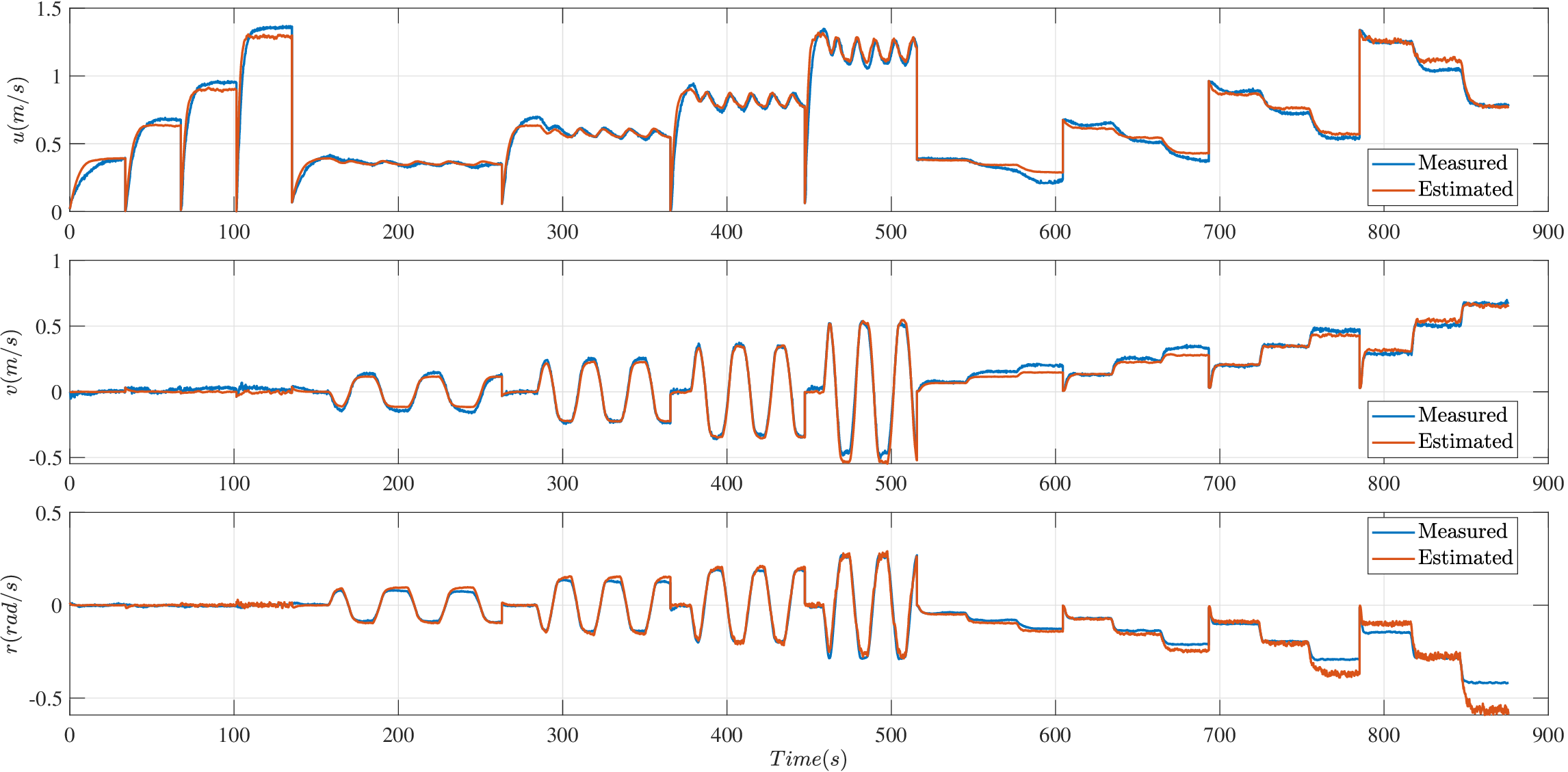}
  \caption{The data validation from the estimated result $P_{GO}$ by using $C_{12}^{12}$ maneuver data.}
  \label{f maneuversvalidation}
\end{figure*}

\begin{table}[htbp]
  \centering
  \caption{Estimated parameter values of Qiuxin No.5 tug.} \label{tbl:The estimated parameters of Qiuxin No.5 tug.}
    \begin{tabular}{cccc}
      \hline
      Parameters & Result & Unit \\
      \hline
      $m_{11}$ & 138.0574  & $kg$ \\
      $m_{22}$ & 106.6003  & $kg$ \\
      $m_{23}$ & 1.1254 & $kg$ \\
      $m_{32}$ & -16.0598  & $kg$ \\
      $m_{33}$ & 15.6476  & $kgm^2$ \\
      $X_u$ & -8.9859   & $kg/s$ \\
      $X_{|u| u}$ & -31.4285  & $kg/s$ \\
      $X_{uuu}$ & -6.8953 & $kg/s$ \\
      $Y_v$ & -71.9041  & $kg/s$ \\
      $Y_{|v|v}$ & -77.6429& $kg/s$ \\
      $Y_{vvv}$ & -27.1394  & $kg/s$ \\
      $Y_{|r|v}$ & -43.2207  & $kg/s$ \\
      $Y_{r}$ & -26.0498  & $kg/s$ \\
      $Y_{|v|r}$ & 26.7652 & $kg/s$ \\
      $Y_{|r|r}$ & 7.7996   & $kg/s$ \\
      $N_v$ & -14.8953  & $kg/s$ \\
      $N_{|v|v}$ & -1.6306   & $kg/s$ \\
      $N_{|r|v}$ & 8.7911   & $kg/s$ \\
      $N_r$ & -26.7122 & $kg/s$ \\
      $N_{|r|r}$ &-9.8284 & $kg/s$ \\
      $N_{rrr}$ & -9.2320  &$kg/s$ \\
      $N_{|v|r}$ & -2.3474  &$kg/s$ \\
      \hline
    \end{tabular}
\end{table}

By using Algorithm~\ref{alg1}, the performance of the estimated parameters compared with the measured data are shown in Figure~\ref{f CA validation} and \ref{f maneuversvalidation}. These figures present the results of the GO and LO processes using $C_{12}^{12}$ maneuver data. The performance of the LO solution is generally consistent with the measured values, except for the $r$ values in the last four turning circle experiments. The primary source of error in the yaw motion can be attributed to the simplification of the transform $\boldsymbol{J}(\psi)$, as only ship states $(x,y,\psi)$ are the optimization target. The results of the GO process exhibit slightly higher errors compared to those of LO, but they remain very close to the measured values. From Figure~\ref{f maneuversvalidation}, it can be observed that the parameters obtained through GO exhibit a consistent motion trend with the actual ship parameters when subjected to identical inputs. Finally, the estimated parameters are shown in Table~\ref{tbl:The estimated parameters of Qiuxin No.5 tug.}.

\subsection{The comparison with the 6-parameter model.}
\begin{figure*}[htbp]
  \centering
  \includegraphics[width = 1\textwidth]{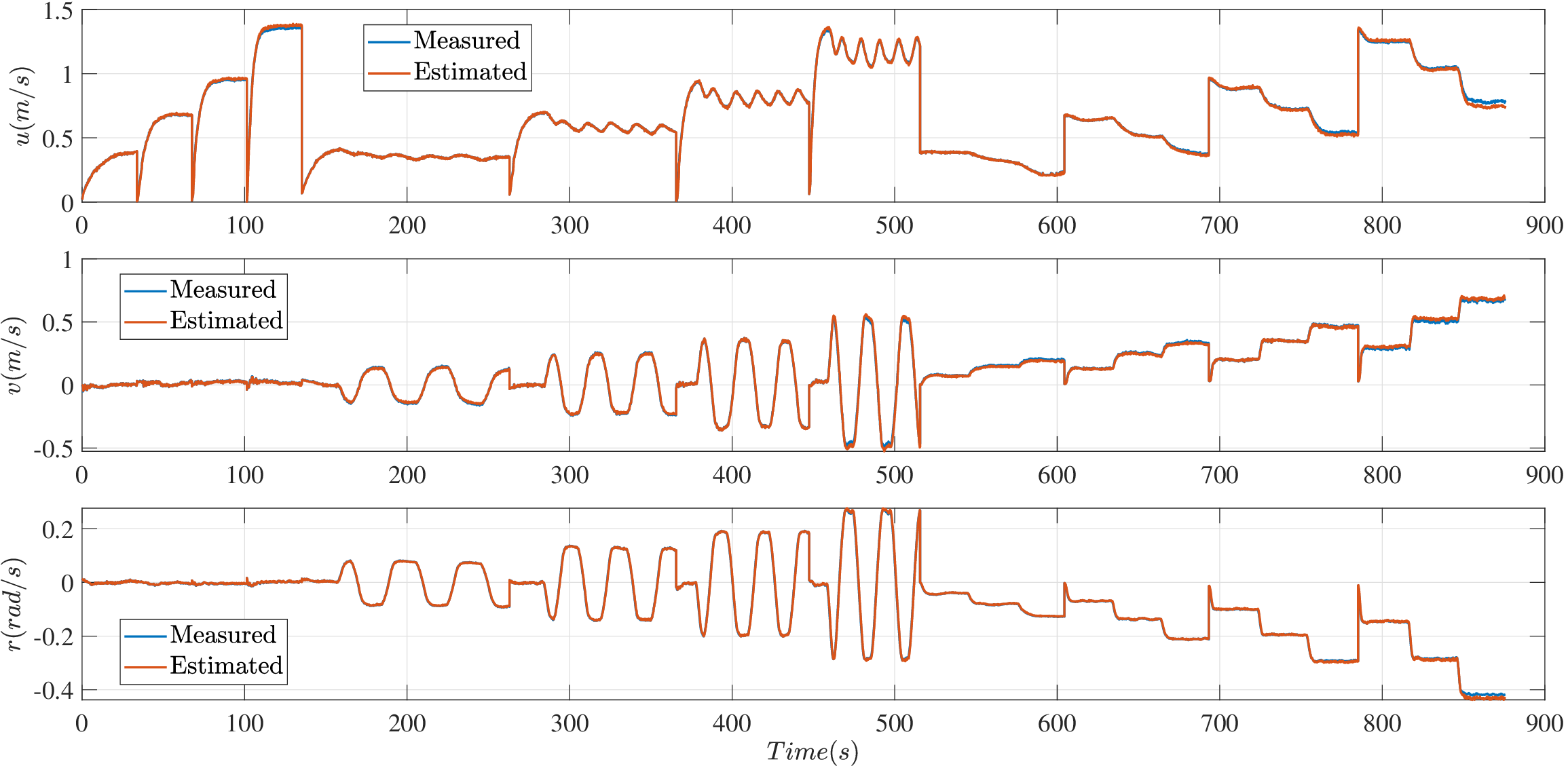}
  \caption{The maneuvers data validation from the estimated result by 6-parameter model with state update in every step.}
  \label{f 6-parameter CA results CA tests}
\end{figure*}

\begin{figure*}[htbp]
  \centering
  \includegraphics[width = 1\textwidth]{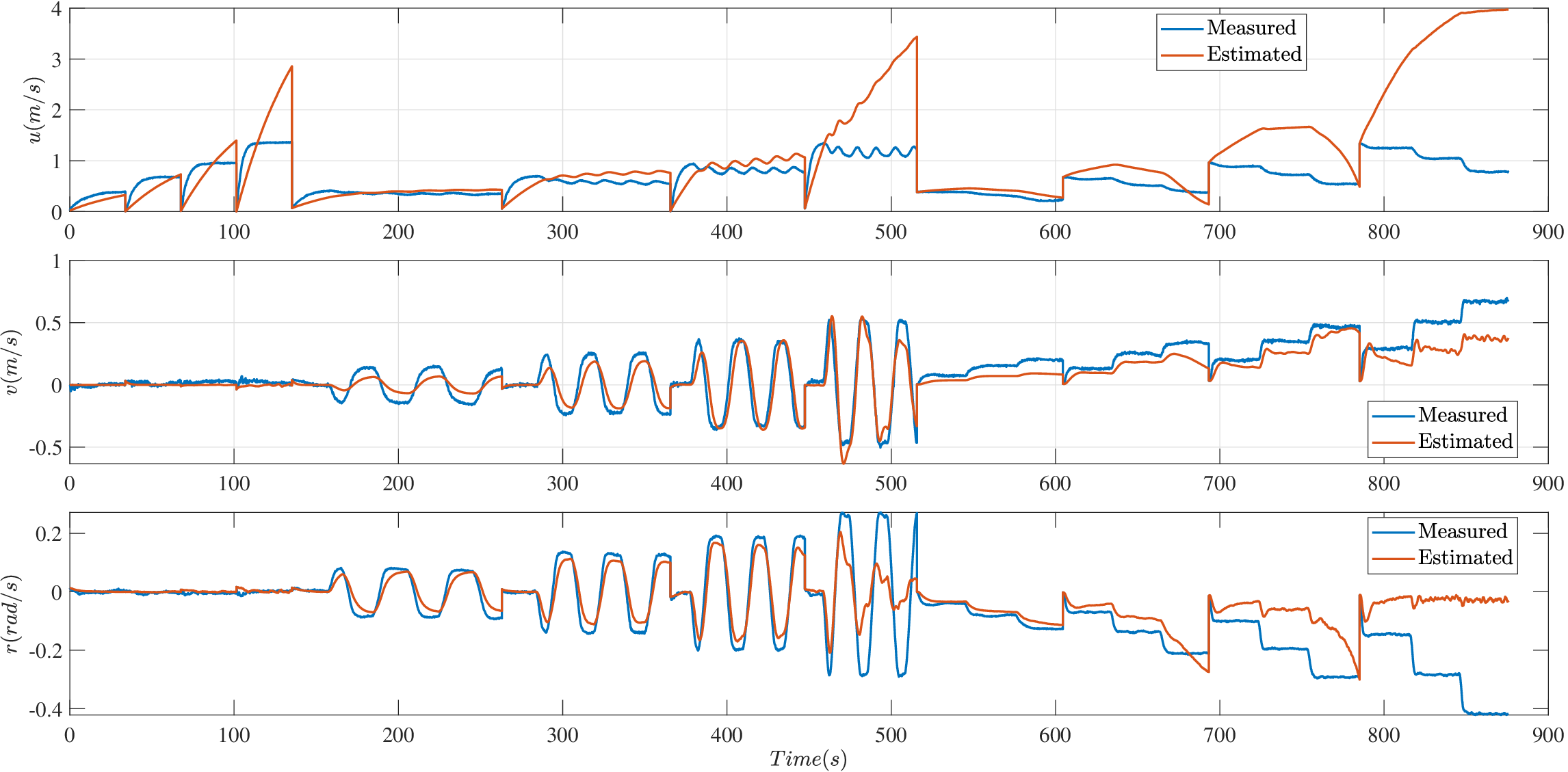}
  \caption{The maneuvers data validation from the estimated result by 6-parameter model without state update.}
  \label{f 6-parameter CA results MBPE tests}
\end{figure*}

Figures~\ref{f 6-parameter CA results CA tests} and \ref{f 6-parameter CA results MBPE tests} illustrate the performance of the estimated parameters using the 6-parameter model,  and the parameter values are shown in Table \ref{tbl:The estimated 6-parameter of Qiuxin No.5 tug.}. It can be observed that the 22-parameter model achieves a better match with the measured data compared to the 6-parameter model in the GO process. Specifically, the 6-parameter model cannot accurately capture the hydrodynamic characteristics of the ship under continuous control input without the state update, leading to larger errors in motion responses according to the comparison between Figure~\ref{f 6-parameter CA results MBPE tests} and Figure~\ref{f maneuversvalidation}. The 6-parameter model fails to accurately simulate the surge movement and yaw movement, causing significant discrepancies in the corresponding response curves. 

\begin{table}[htbp]
  \centering
  \caption{Estimated parameter values of Qiuxin No.5 tug for the 6-parameter model by using $C_{12}^{12}$ maneuver data.} \label{tbl:The estimated 6-parameter of Qiuxin No.5 tug.}
    \begin{tabular}{cccc}
      \hline
      Parameters & Result & Unit \\
      \hline
      $m_{11}$ & 216.4727  & $kg$ \\
      $m_{22}$ & 183.4906  & $kg$ \\
      $m_{33}$ & 0.0632 & $kgm^2$ \\
      $d_{11}$ & 44.1878  & $kg/s$ \\
      $d_{22}$ & 152.9380 & $kg/s$ \\
      $d_{33}$ & 0.2629 & $kg/s$ \\
      \hline
    \end{tabular}
\end{table}

Table~\ref{tbl:relativeError} quantitatively presents the relative errors between the measured data and the simulated results using the 22-parameter model and the 6-parameter model. It is evident that the 6-parameter model results in significantly larger relative errors compared to the 22-parameter model for most of the variables, indicating its limited accuracy in capturing the complete dynamics of the ship. For example, the relative errors in surge velocity, sway velocity, and yaw velocity for the 6-parameter model are 215\%, 344\%, and 445\%, respectively, while the corresponding values for the 22-parameter model are only 2.33\%, 4.45\%, and 5.85\%, respectively. These results further confirm the superiority of the 23-parameter model over the 6-parameter model in accurately representing the ship maneuvering behavior.
\begin{table}[htbp]
  \centering
  \caption{Relative errors between measured data and simulated results.} \label{tbl:relativeError}
    \begin{tabular}{cccccc}
      \hline
      Variables & 22-Parameter Model & 6-Parameter Model \\
      \hline
      Surge Vel. & 2.33\% & 215.4\% \\
      Sway Vel. & 4.45\% & 344.5\% \\
      Yaw Vel. & 5.85\% & 445.5\% \\
      Surge Acc. & 0.83\% & 46.2\% \\
      Sway Acc. & 2.52\% & 85.6\% \\
      Yaw Acc. & 2.80\% & 89.9\% \\
      Surge Dist. & 7.16\% & 110.3\% \\
      Sway Dist. & 17.4\% & 225.7\% \\
      Yaw Dist. & 15.9\% & 329.8\% \\
      \hline
    \end{tabular}
\end{table}

Figure \ref{f 6-parameter CA results CA tests} demonstrates that the predictive accuracy of the six-parameter model can be achieved with state updates at each step. This implies that the 6-parameter model is only suitable for short iteration prediction rather than long iteration prediction. Furthermore, it indicates that the 22-parameter model aligns more closely with the true motion process of the tug. In conclusion, these results underscore the importance of considering a more comprehensive model when estimating ship dynamic parameters.




\subsection{Thruster model identification}
The thruster model identification involves considering the speed acceleration model, the thruster force model, and the thrust configuration model. As the target ship exhibits a very rapid speed response between the command and feedback, the speed acceleration model is omitted. The other two models are described below.

\subsubsection{ Thruster force regression}
\begin{figure}[htbp]
  \centering
  \includegraphics[width = 0.45\textwidth]{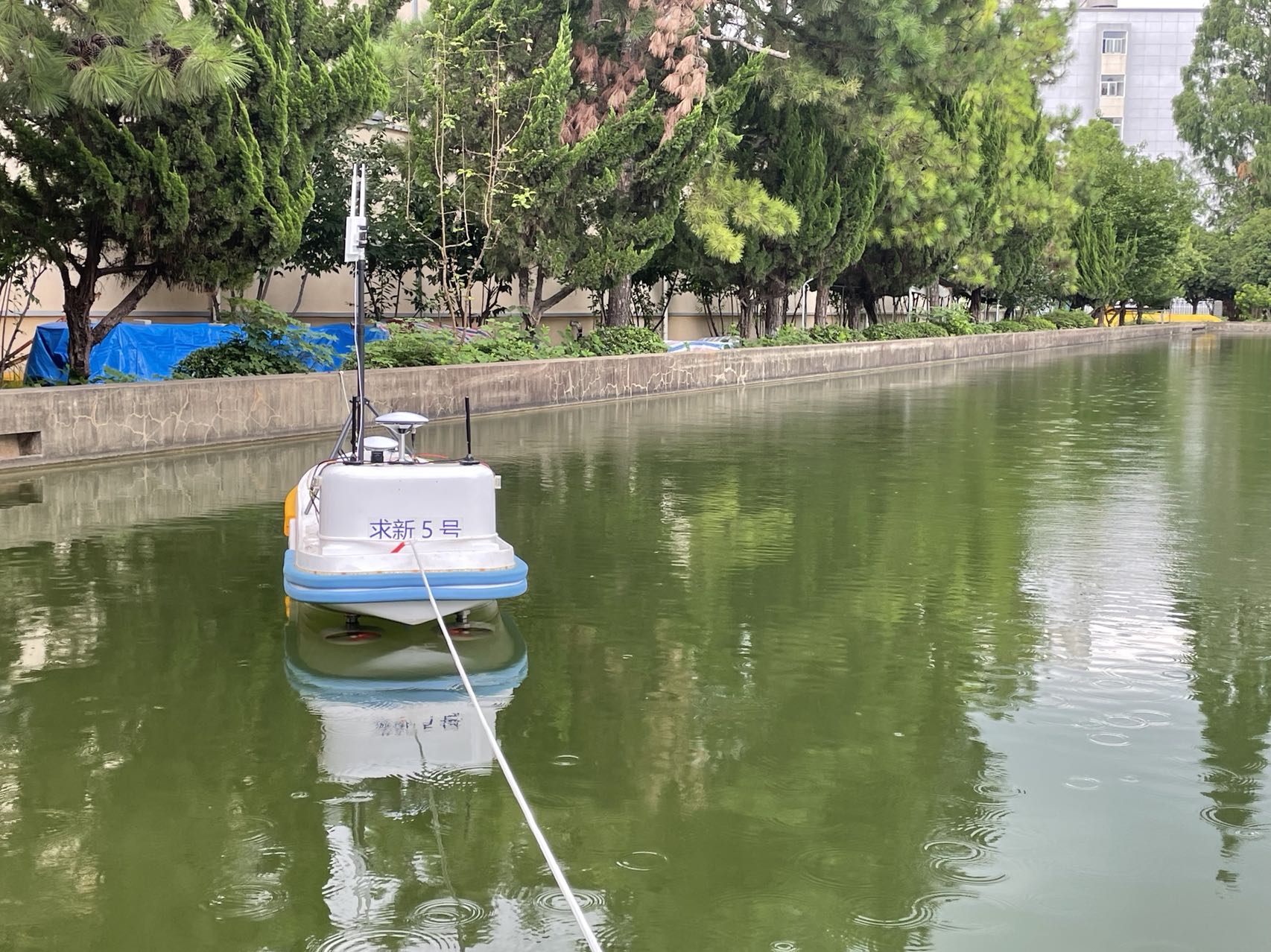}
  \caption{The simplified bollard pull test of the tug.}
  \label{f bollard pull test}
\end{figure}

\begin{figure}[htbp]
  \centering
  \includegraphics[width = 0.45\textwidth]{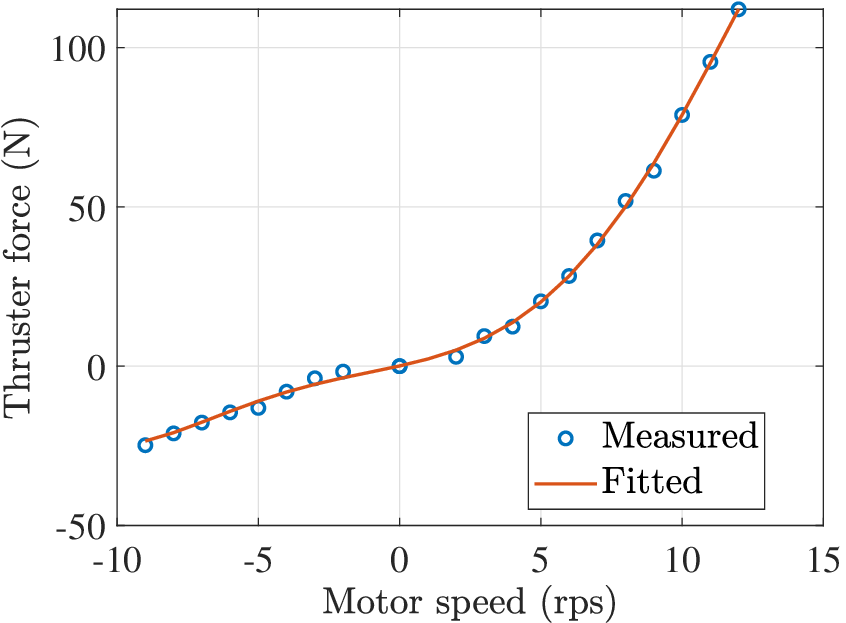}
  \caption{Curve fitting between RPS and thrusters.}
  \label{f:curve fitting between RPS and force}
\end{figure}

According to the thruster configuration mentioned in Section \ref{sec2}, the relation between thruster motor speed and thruster force is obtained through the bollard pull test, as shown in Figure \ref{f bollard pull test}. A 5th-order polynomial regression is performed to fit the data, giving the following equation:
\begin{equation} \label{eq20}
  \begin{aligned}
    f=& (-0.0001773n^5 + 0.001187n^4 \\
    &+ 0.04978n^3 + 0.151n^2 +1.974n + 0.0722)/2
  \end{aligned}
\end{equation}
where $n$ represents the thruster speed, and $f$ represents the calculated force. The results are shown in Figure \ref{f:curve fitting between RPS and force}.

\subsubsection{Azimuth angular acceleration model identification}

\begin{figure*}[htbp]
  \centering
  \includegraphics[width = 1\textwidth]{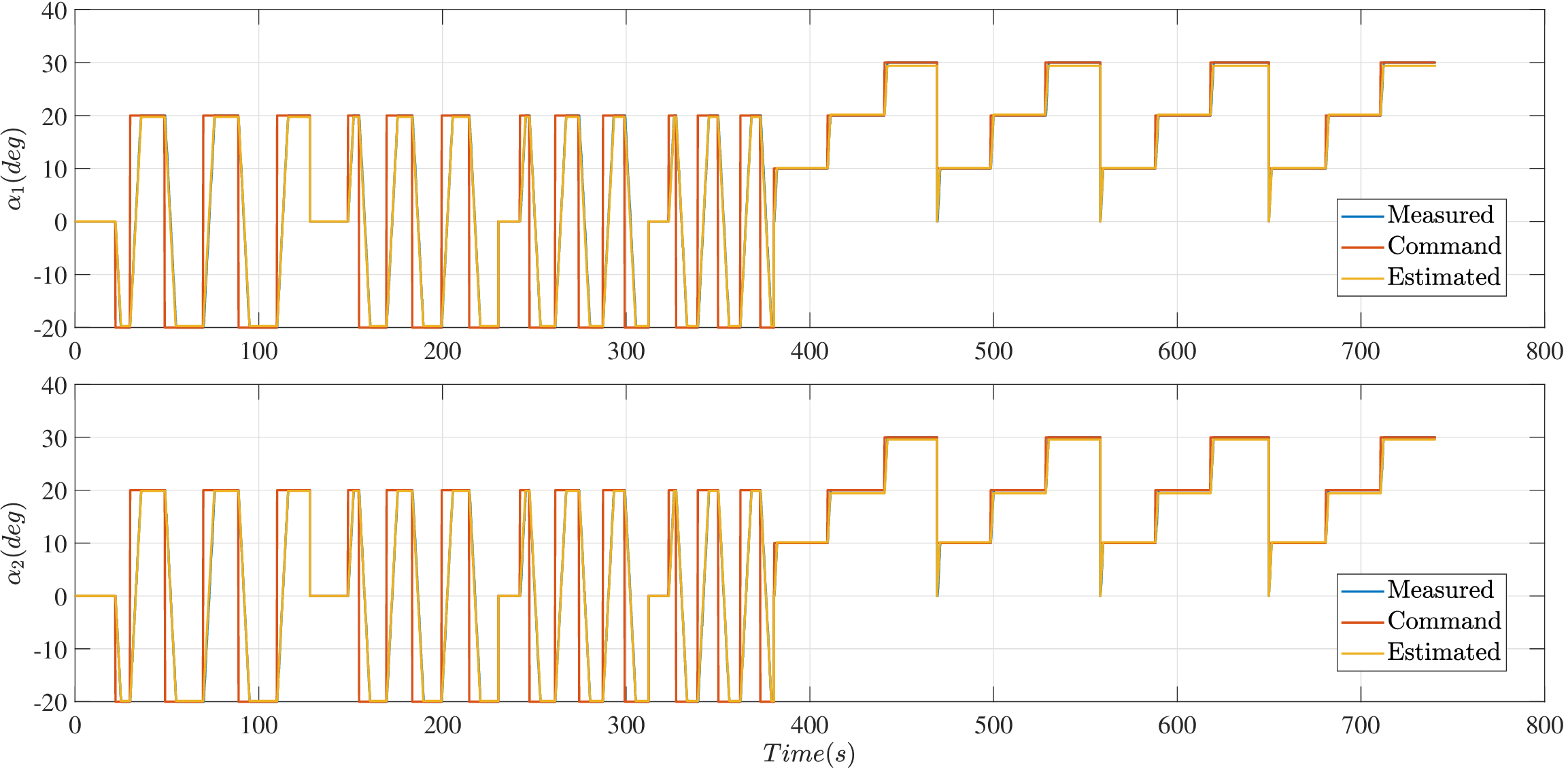}
  \caption{Experimental data and estimated azimuth angles after the parameters have been identified for the azimuth-angle model.}
  \label{f azimuth angel model identification}
\end{figure*}

The model referred to in Section \ref{sec2} is identified using the 5 to 12 maneuvers as test data. By using the GO method, the parameter of the azimuth thruster rotation angle model is obtained. Due to $\epsilon_1 = 0$ and $ \epsilon_2 = 0$, the angular velocity can thus be represented as $\dot{\alpha} = K_{\alpha} \text{sign}(\alpha_d - \alpha) $, where $ K_1 = 0.1151 $ and $ K_2 = 0.1161 $. The performance of the identified model is demonstrated in Figure \ref{f azimuth angel model identification}, where the identified model fits the measured data well. The root mean square error (RMSE) between $\alpha_1$ and $\alpha_2$ is only 0.9992 and 0.9994, respectively, indicating that the identified model has similar response characteristics as the original azimuth angle. 

By using the thruster formula to calculate the thruster force and the model estimation method to determine the rotation angle acceleration model, the final input torques can be calculated based on the thruster configuration matrix given by Equations \ref{eq6} and \ref{eq7}. The parameters in Equation \ref{eq7} are determined by the position of the two azimuth thrusters, as follows: $lx_1 = -0.8$ m, $ lx_2 = -0.8$ m, $ ly_1 = 0.163$ m, and $ ly_2 = -0.163$ m.


In conclusion, the structure of the proposed whole ship model, from input command to output earth-fixed position, is demonstrated in Figure \ref{f the whole ship motion model}. It can be divided into two parts: the input torque model, which includes the calculation of motor speed and force, the model of rotation angle acceleration, and the thruster configuration matrix; and the ship dynamic model, which includes the kinematic model and kinetic model. Therefore, the whole ship motion model can be summarized as:
\begin{equation}
  \dot{\boldsymbol{x}} = f(\boldsymbol{x},\boldsymbol{u}),
\end{equation}
where $\boldsymbol{x} = [\alpha_1,\alpha_2,x,y,\psi,u,v,r]^{\rm{T}}\in \Re^8 $ denotes the ship state vector, and $\boldsymbol{u} =[n_1,n_2,\alpha_1,\alpha_2]^{\rm{T}}\in \Re^4 $ represents the actual input command, in which $n_1,n_2$ refer to the propellers revolution rate and $\alpha_1,\alpha_2$ refer to the azimuth rotation angles. The intermediate control variables $f_1,f_2,\boldsymbol{\tau} $ are calculated using Equation \ref{eq20}, \ref{eq6} and \ref{eq7}.

\section{Model validation results and discussions} \label{sec5}
In this section, we validate the ship motion model using experimental data collected from a real ship. The validation consists of two parts. In the first part, we use estimated data and related zigzag data with the inverse initial heading angle to validate the accuracy of the data under one-step prediction. This part reflects the accuracy of the estimated model under the LO method (1-step prediction). The second part investigates the prediction under different step sizes to verify the effectiveness of the GO method (5-step and longer-step prediction).

\subsection{ Estimation data validation}
\begin{figure*}[htbp]
  \subfloat[+10+10 zigzag]{\includegraphics[width=0.3\textwidth]{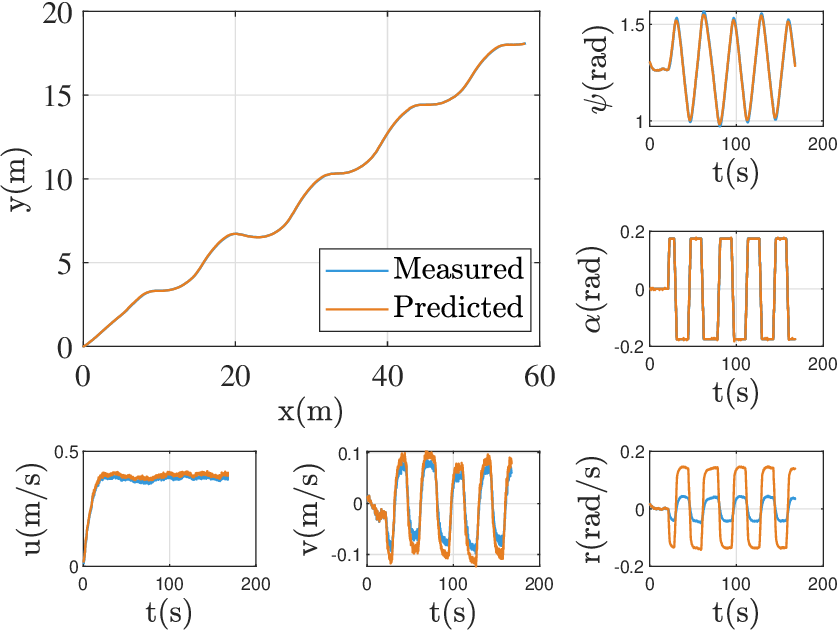}}\label{f RPS3_zigzag1010}
 \hfill 	
  \subfloat[+20+20 zigzag]{\includegraphics[width=0.3\textwidth]{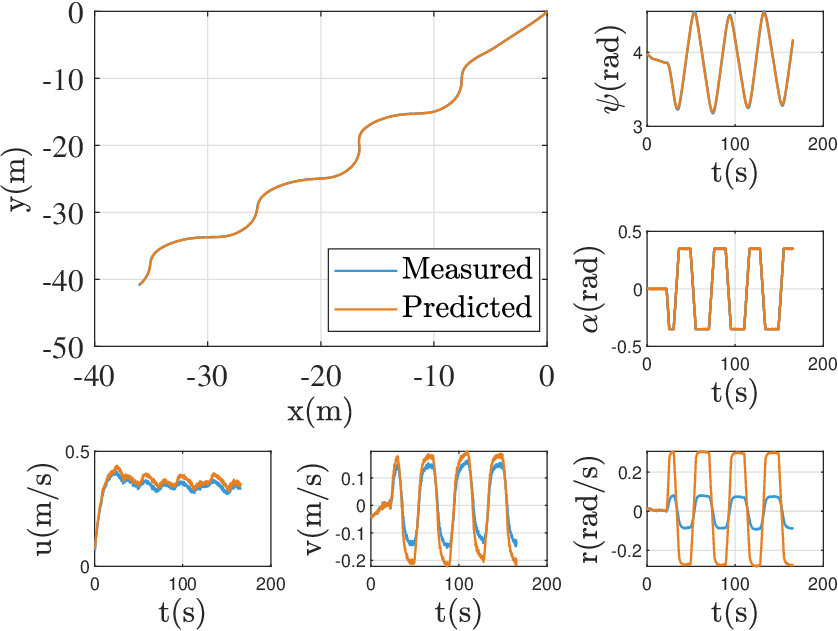}}\label{f RPS3_zigzag2020}
 \hfill	
  \subfloat[+30+30 zigzag]{\includegraphics[width=0.3\textwidth]{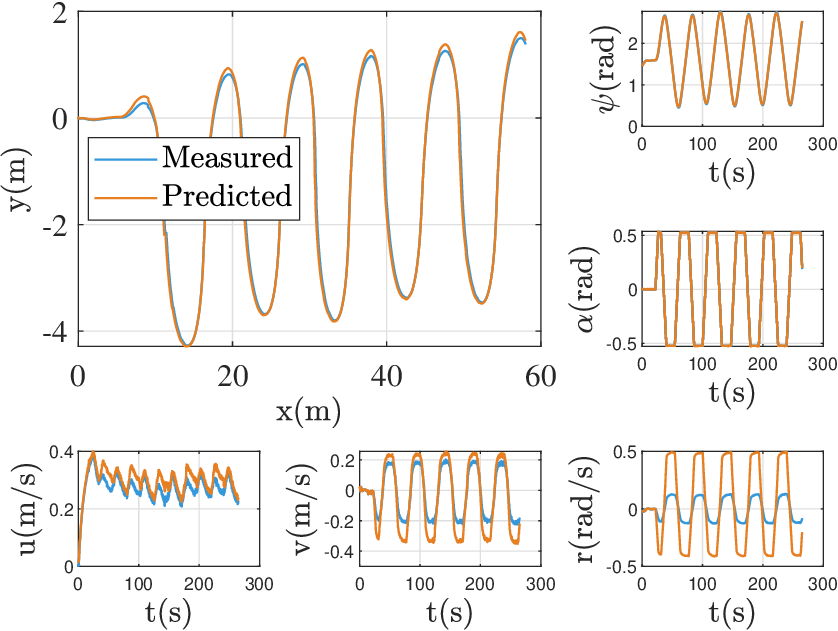}}\label{f RPS3_zigzag3030}
  \newline
  \subfloat[+10-10 zigzag]{\includegraphics[width=0.3\textwidth]{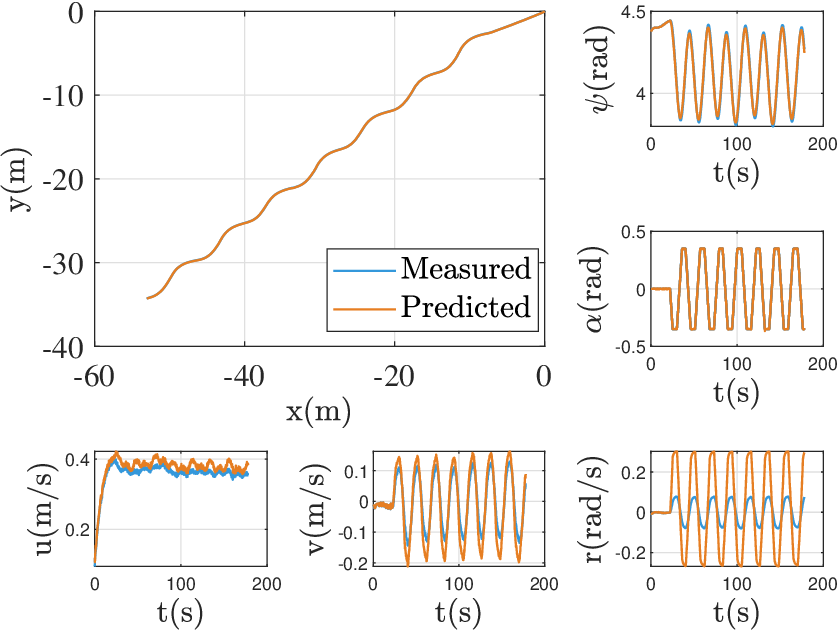}}
  \label{f RPS3_zigzag10_10}
 \hfill 	
  \subfloat[+20-20 zigzag]{\includegraphics[width=0.3\textwidth]{figure/05prediction02/RPS3_zigzag20_20.mat.eps}}
  \label{f RPS3_zigzag20_20}
 \hfill	
  \subfloat[+30-30 zigzag]{\includegraphics[width=0.3\textwidth]{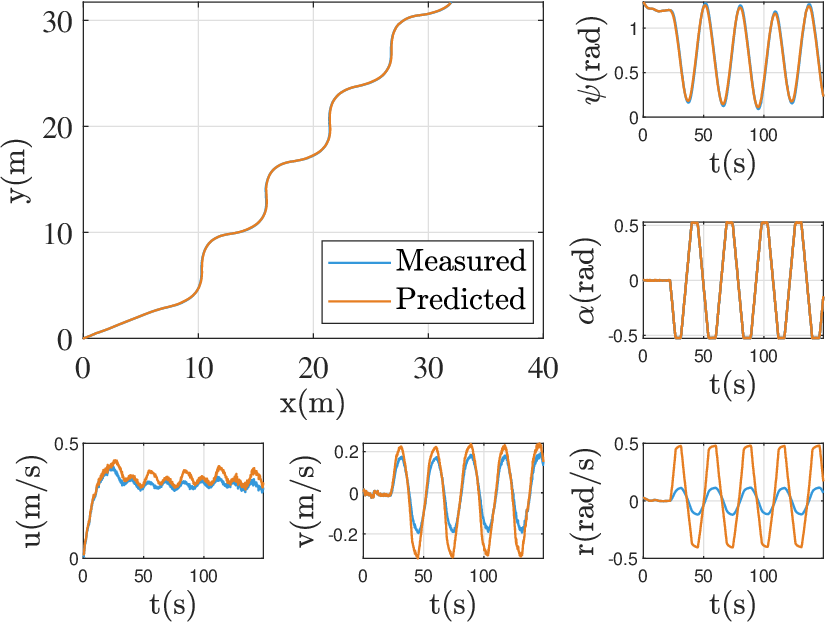}}
  \label{f RPS3_zigzag30_30}
\caption{The validation based on 3 RPS motor speed.}
\end{figure*}
\begin{figure*}[htbp]
  \subfloat[+10+10 zigzag]{\includegraphics[width=0.3\textwidth]{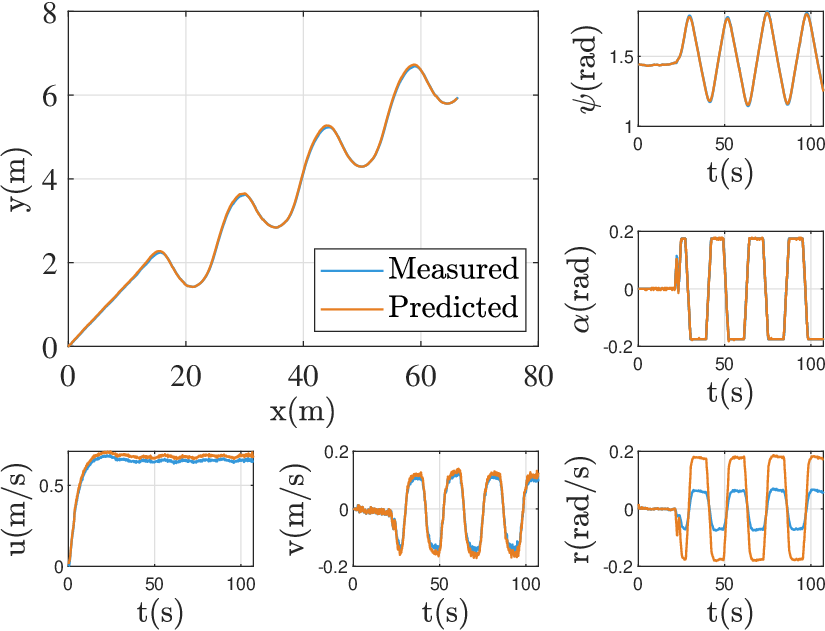}}
 \hfill 	
  \subfloat[+20+20 zigzag]{\includegraphics[width=0.3\textwidth]{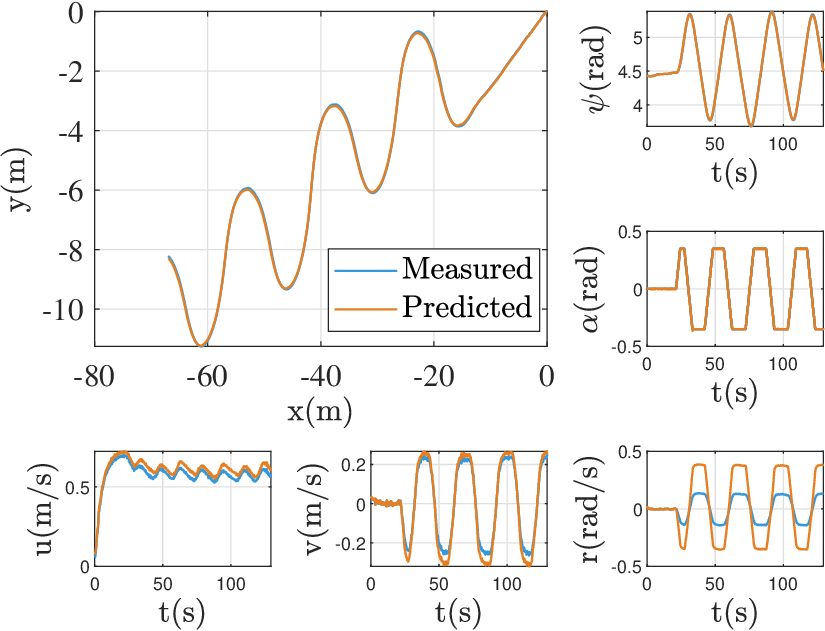}}
 \hfill	
  \subfloat[+30+30 zigzag]{\includegraphics[width=0.3\textwidth]{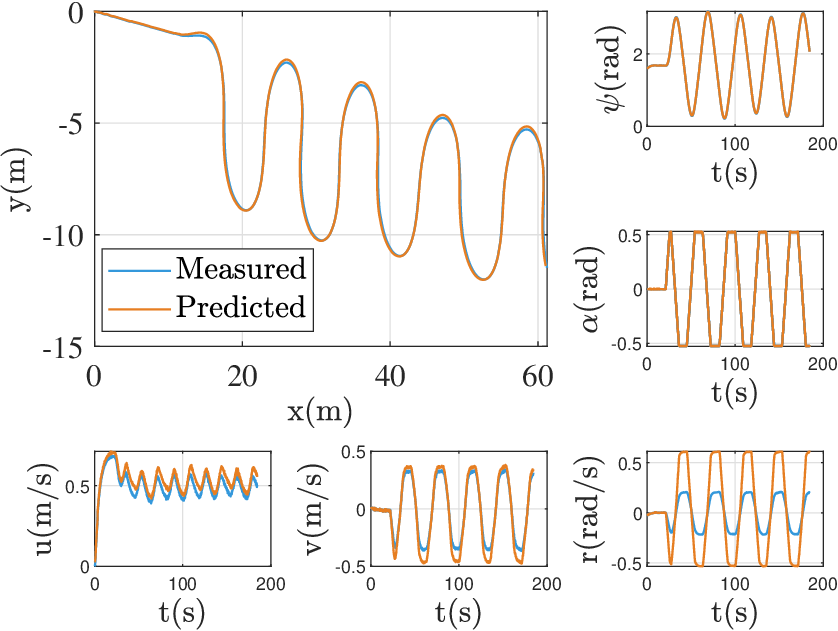}}
  \newline
  \subfloat[+10-10 zigzag]{\includegraphics[width=0.3\textwidth]{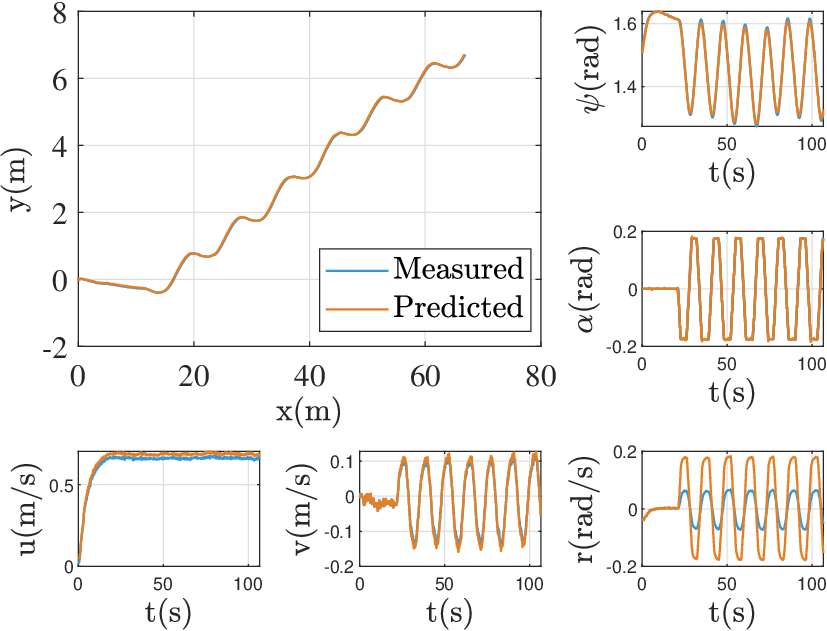}}
  \label{f RPS5_zigzag10_10}
 \hfill 	
  \subfloat[+20-20 zigzag]{\includegraphics[width=0.3\textwidth]{figure/05prediction02/RPS5_zigzag10_10.mat.eps}}
  \label{f RPS5_zigzag20_20}
 \hfill	
  \subfloat[+30-30 zigzag]{\includegraphics[width=0.3\textwidth]{figure/05prediction02/RPS5_zigzag10_10.mat.eps}}
  \label{f RPS5_zigzag30_30}
\caption{The validation based on 5 RPS motor speed.}
\end{figure*}
\begin{figure*}[htbp]
  \subfloat[+10+10 zigzag]{\includegraphics[width=0.3\textwidth]{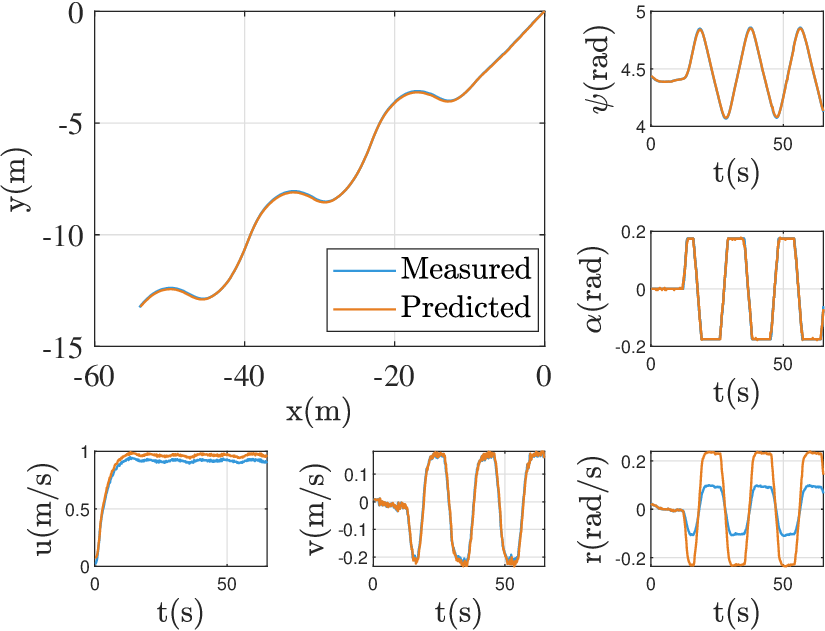}}\label{f RPS7_zigzag1010}
 \hfill 	
  \subfloat[+20+20 zigzag]{\includegraphics[width=0.3\textwidth]{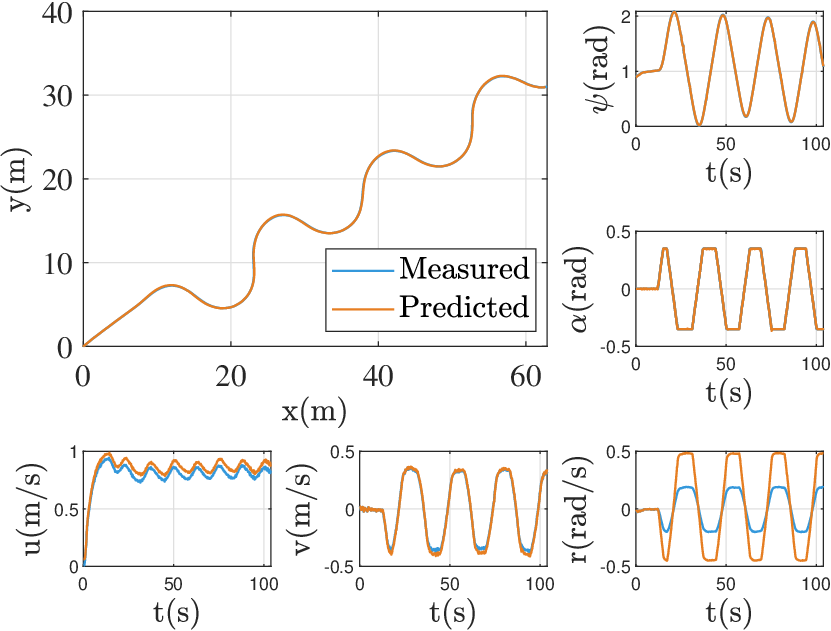}}\label{f RPS7_zigzag2020}
 \hfill	
  \subfloat[+30+30 zigzag]{\includegraphics[width=0.3\textwidth]{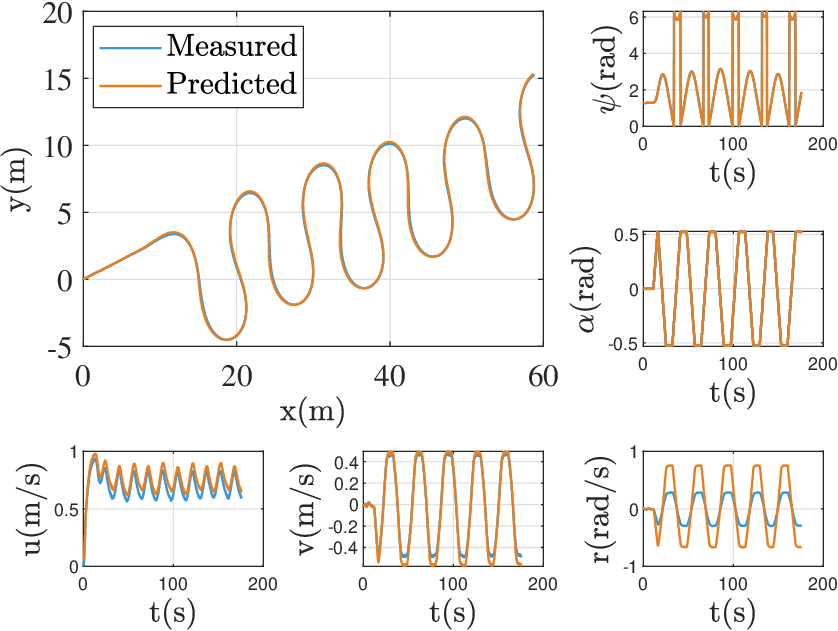}}\label{f RPS7_zigzag2020}
  \newline
  \subfloat[+10-10 zigzag]{\includegraphics[width=0.3\textwidth]{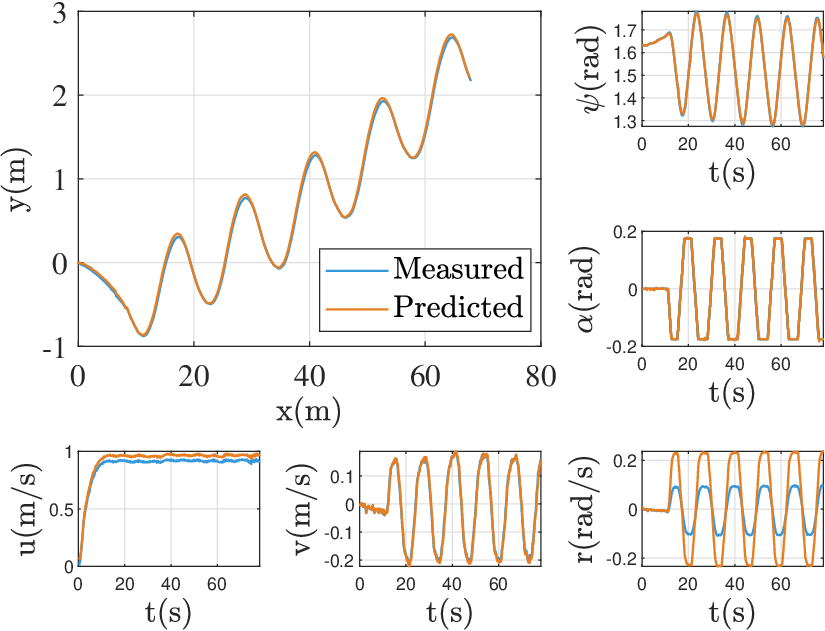}}
  \label{f RPS7_zigzag10_10}
 \hfill 	
  \subfloat[+20-20 zigzag]{\includegraphics[width=0.3\textwidth]{figure/05prediction02/RPS7_zigzag10_10.mat.eps}}
  \label{f RPS7_zigzag20_20}
 \hfill	
  \subfloat[+30-30 zigzag]{\includegraphics[width=0.3\textwidth]{figure/05prediction02/RPS7_zigzag10_10.mat.eps}}
  \label{f RPS7_zigzag30_30}
\caption{The validation based on 7 RPS motor speed.}
\end{figure*}
\begin{figure*}[htbp]
  \subfloat[+10+10 zigzag]{\includegraphics[width=0.3\textwidth]{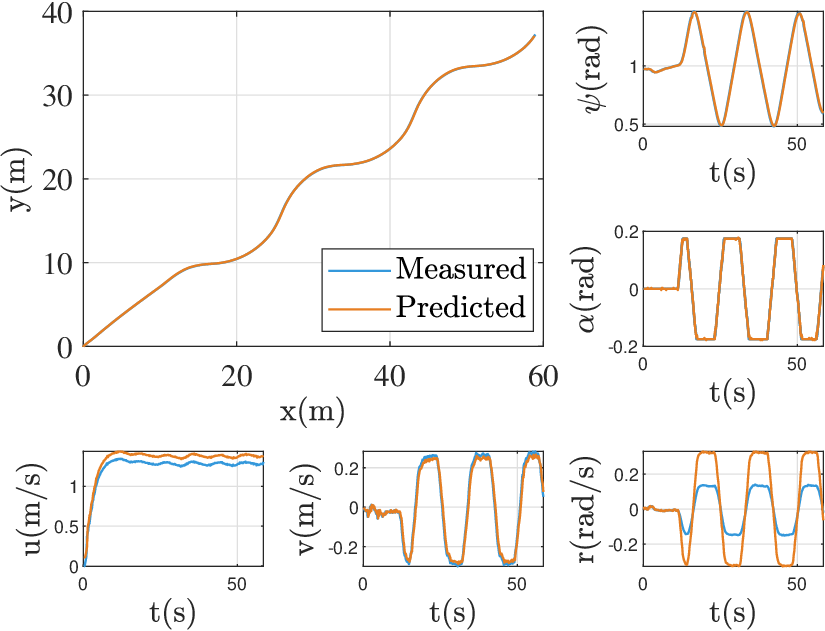}}\label{f RPS10_zigzag1010}
 \hfill 	
  \subfloat[+20+20 zigzag]{\includegraphics[width=0.3\textwidth]{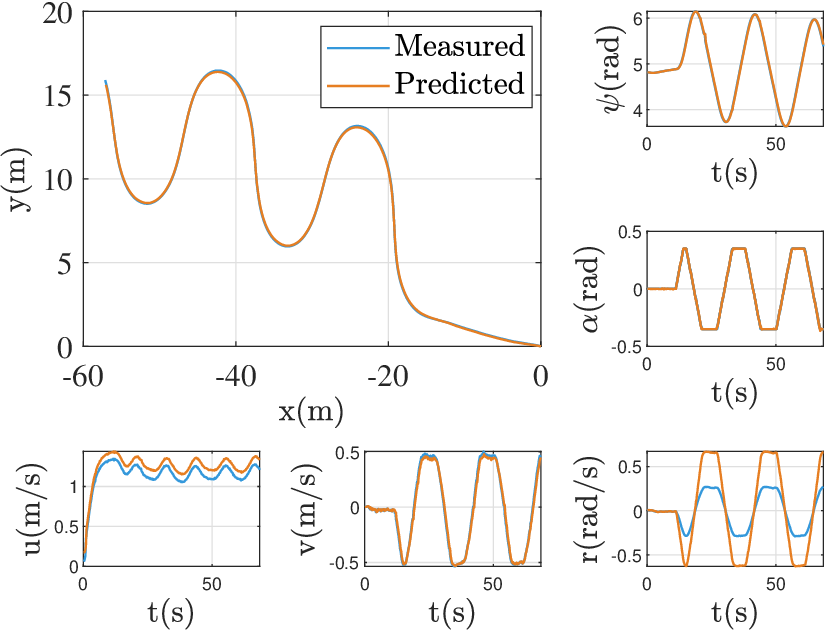}}\label{f RPS10_zigzag2020}
 \hfill	
  \subfloat[+30+30 zigzag]{\includegraphics[width=0.3\textwidth]{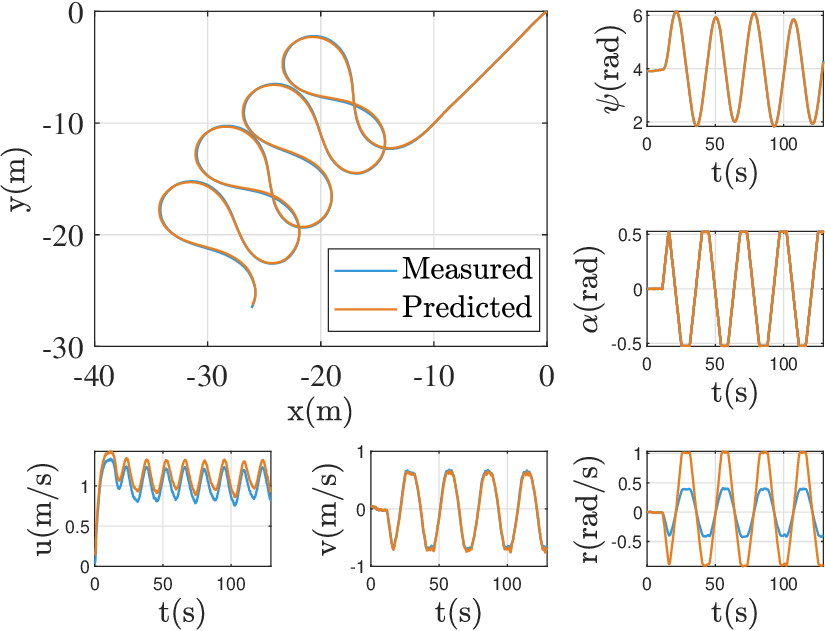}}\label{f RPS10_zigzag2020}
  \newline
  \subfloat[+10-10 zigzag]{\includegraphics[width=0.3\textwidth]{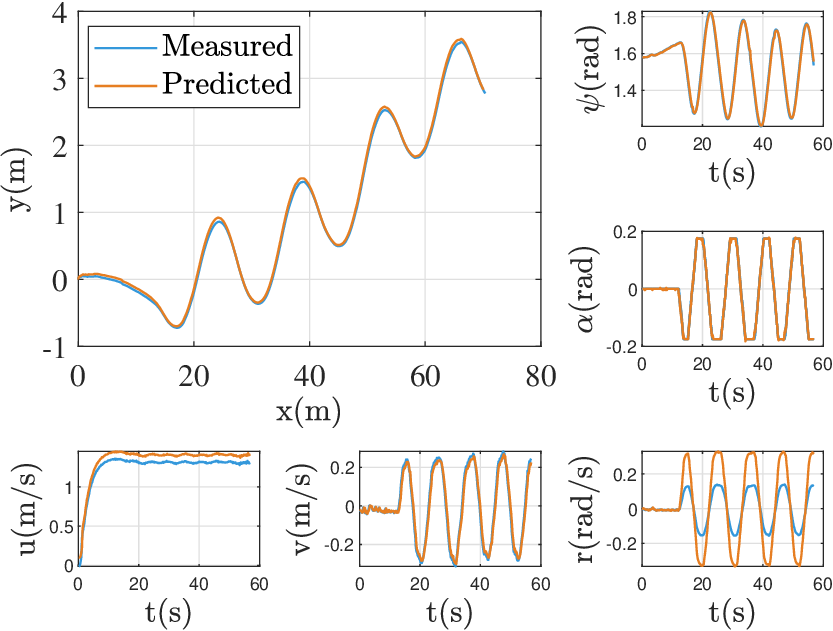}}
  \label{f RPS10_zigzag10_10}
 \hfill 	
  \subfloat[+20-20 zigzag]{\includegraphics[width=0.3\textwidth]{figure/05prediction02/RPS10_zigzag10_10.mat.eps}}
  \label{f RPS10_zigzag20_20}
 \hfill	
  \subfloat[+30-30 zigzag]{\includegraphics[width=0.3\textwidth]{figure/05prediction02/RPS10_zigzag10_10.mat.eps}}
  \label{f RPS10_zigzag30_30}
\caption{The validation based on 10 RPS motor speed.}
\end{figure*}

\begin{figure*}[htbp]
  \subfloat[3 RPS ]{\includegraphics[width=0.24\textwidth]{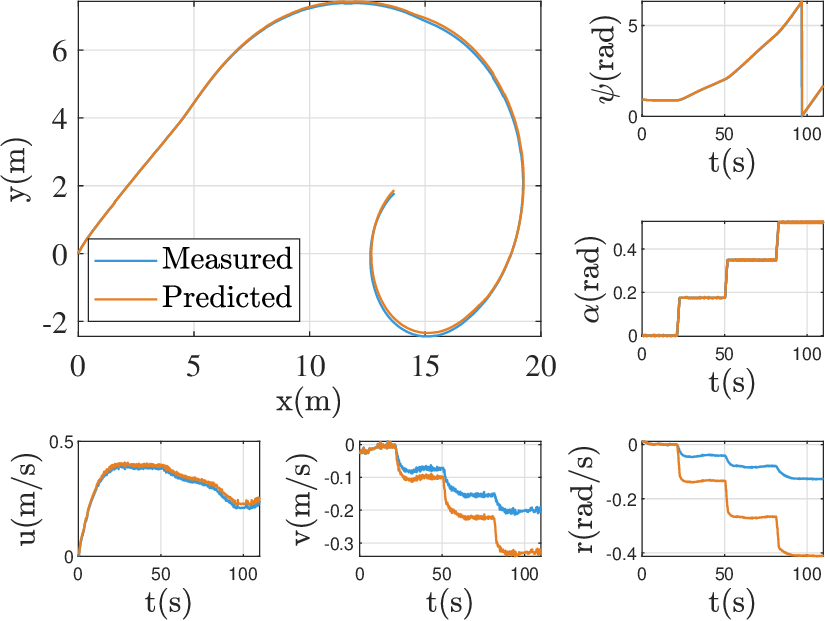}}\label{f RPS3_turningCircle102030 }
 \hfill 	
  \subfloat[5 RPS]{\includegraphics[width=0.24\textwidth]{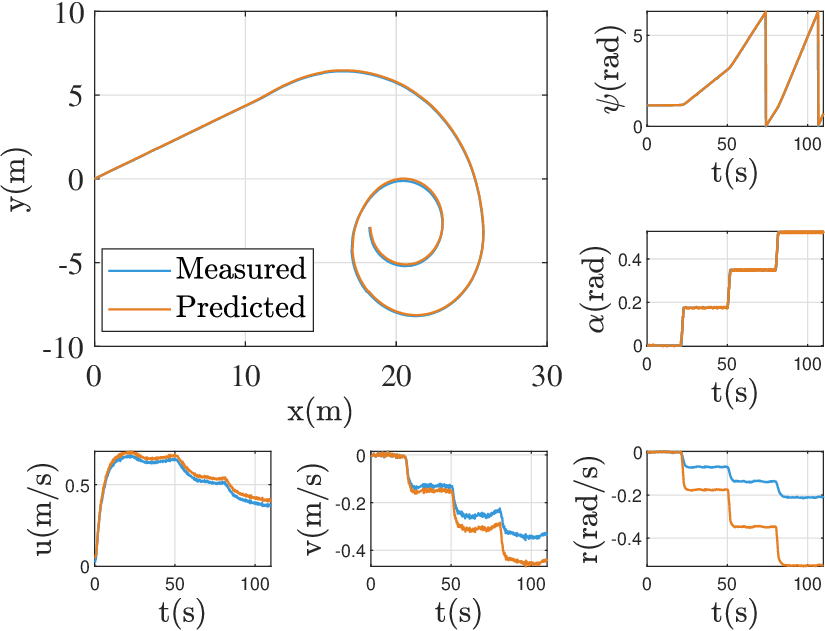}}\label{f RPS5_turningCircle102030}
 \hfill	
  \subfloat[7 RPS]{\includegraphics[width=0.24\textwidth]{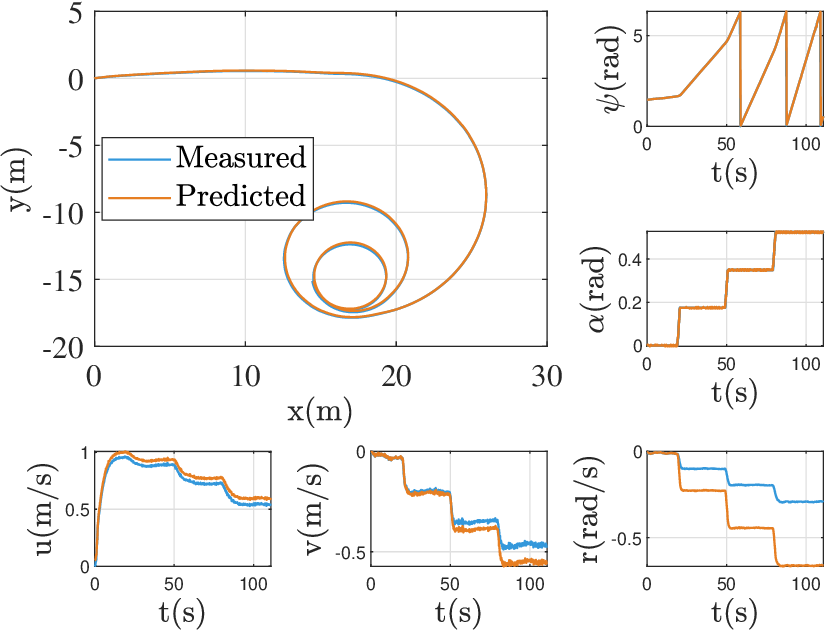}}\label{f RPS7_turningCircle102030}
  \hfill
  \subfloat[10 RPS]{\includegraphics[width=0.24\textwidth]{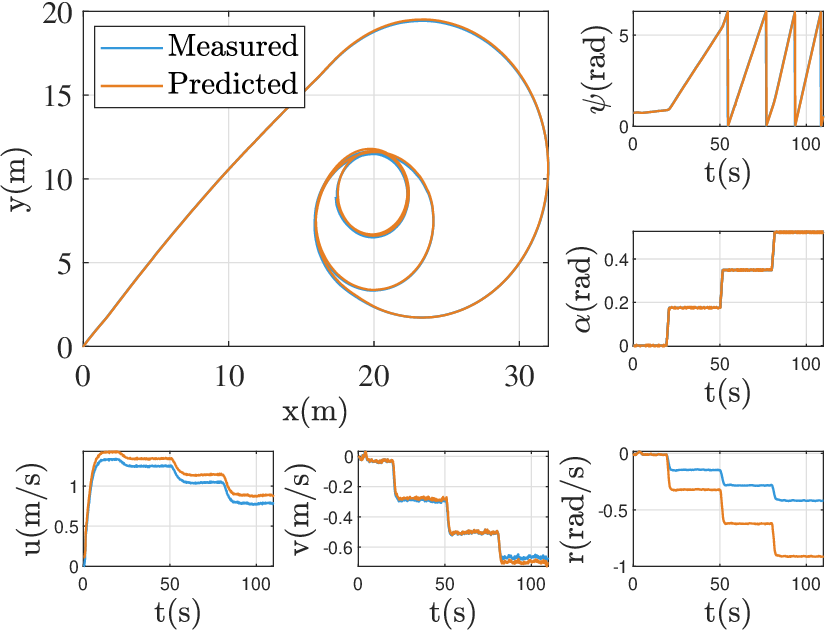}}
  \label{f RPS10_turningCircle102030}
\caption{The validation based on turning circle.}
\end{figure*}
For the validation of the estimated data, we selected straight line movements, zigzag tests, and turning tests. These data were collected from a real tug. Using the RK4 method, we predicted the next states given the states of each step and input. Figures \ref{f RPS3_zigzag1010} to \ref{f RPS10_turningCircle102030} in Appendix \ref{appendix} show the comparisons between the predicted positions $[x,y]$, heading angle $\psi $, and azimuth angles $\alpha$ with the measured real data. The surge and sway speeds are calculated very closely, while the predicted yaw motion has a larger error but follows the same trend.

\begin{figure*}[htbp]
  \centering
  \includegraphics[width = 1\textwidth]{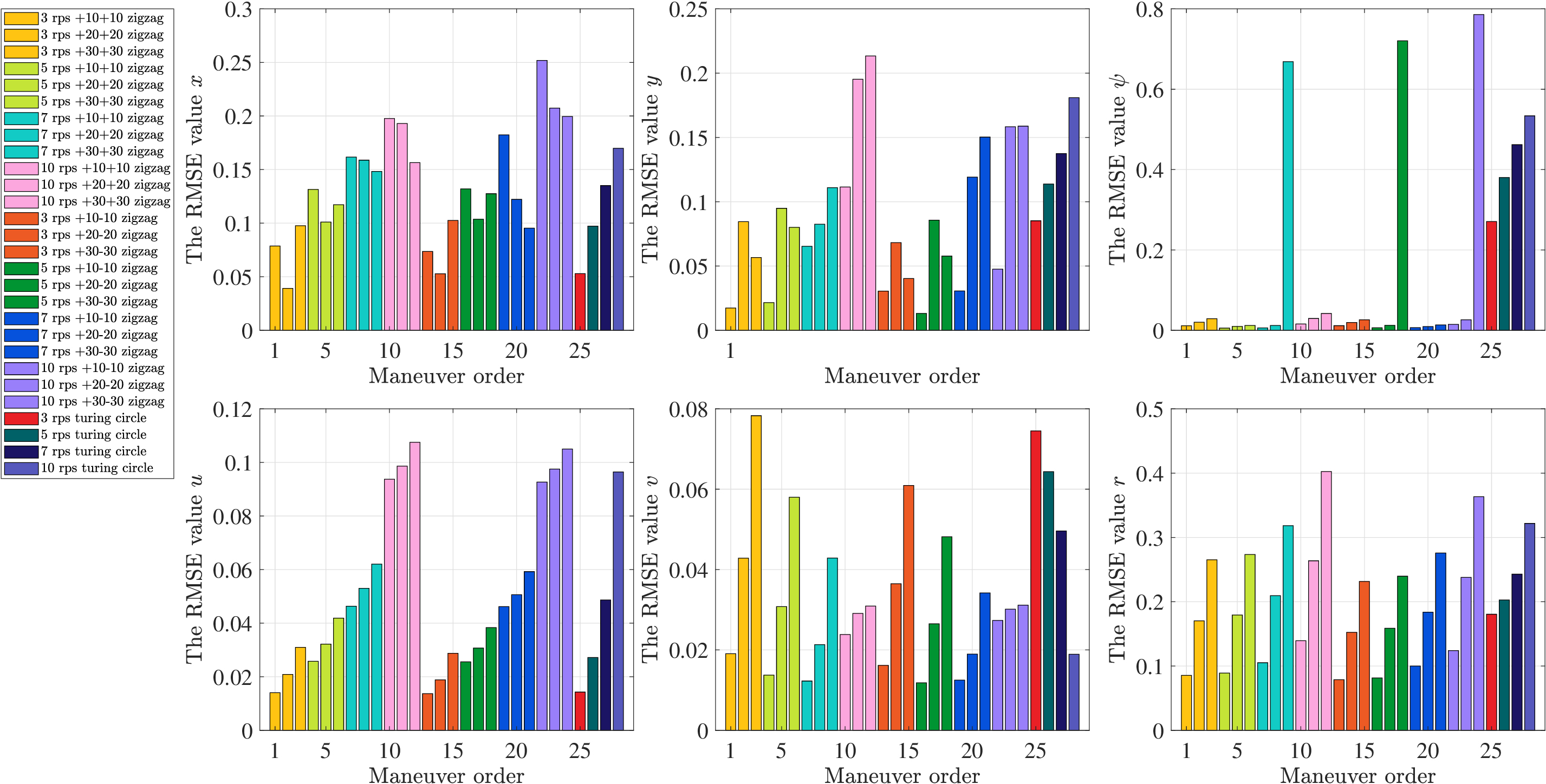}
  \caption{The RMSE values of the 28 maneuvers.}
  \label{f RMSE of the maneuvers.}
\end{figure*}

Figure \ref{f RMSE of the maneuvers.} presents the Root Mean Squared Error (RMSE) values of the total 28 maneuvers, ordered according to Figures \ref{f RPS3_zigzag1010} to \ref{f RPS10_turningCircle102030}. Generally, the RMSE increases with the speed and rotation angle. The error in yaw motion ($r$) is greater than in the other two motions. This difference is likely due to the simplification method used in converting from the Earth fixed coordinate system to the BODY fixed coordinate system, which reduces the accuracy of the yaw velocity prediction while maintaining the accuracy of the Earth fixed coordinate system. The presence of crab angles during the turning circle validation also leads to greater prediction errors in sway and yaw speeds, with faster speeds resulting in larger errors. In conclusion, the experimental results show that the proposed model accurately predicts the ship's position, heading angle, and azimuth angles. The validation demonstrates the effectiveness of the proposed model in ship motion prediction. 

\begin{figure*}[htbp]
  \centering
  \includegraphics[width = 1\textwidth]{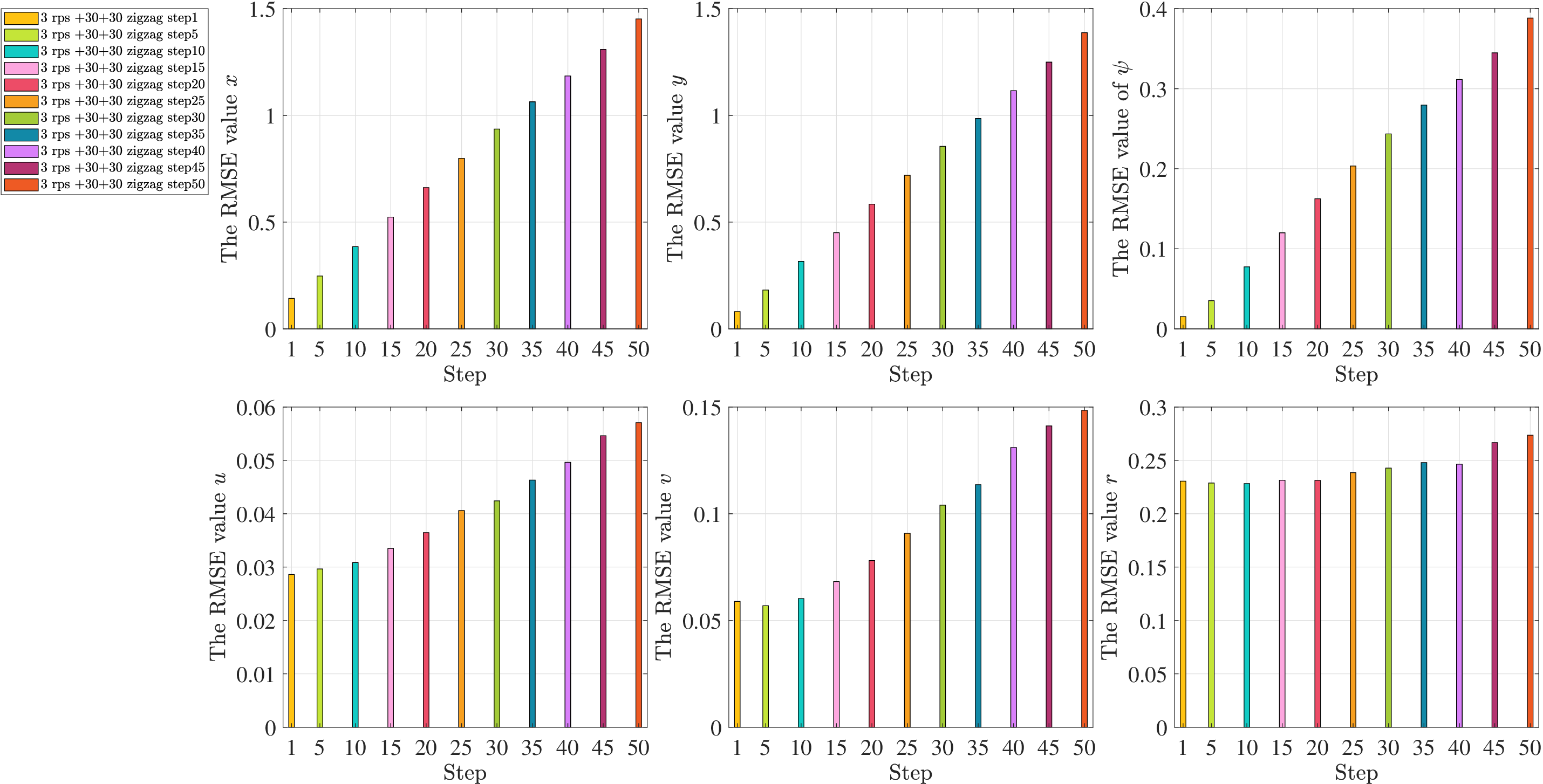}
  \caption{Comparison of RMSE values with different prediction steps for the same measurement data.}
  \label{f RMSE of different steps}
\end{figure*}


\subsection{Model movement prediction validation}
\begin{figure*}[htbp]
  \subfloat[1 step]{\includegraphics[width=0.3\textwidth]{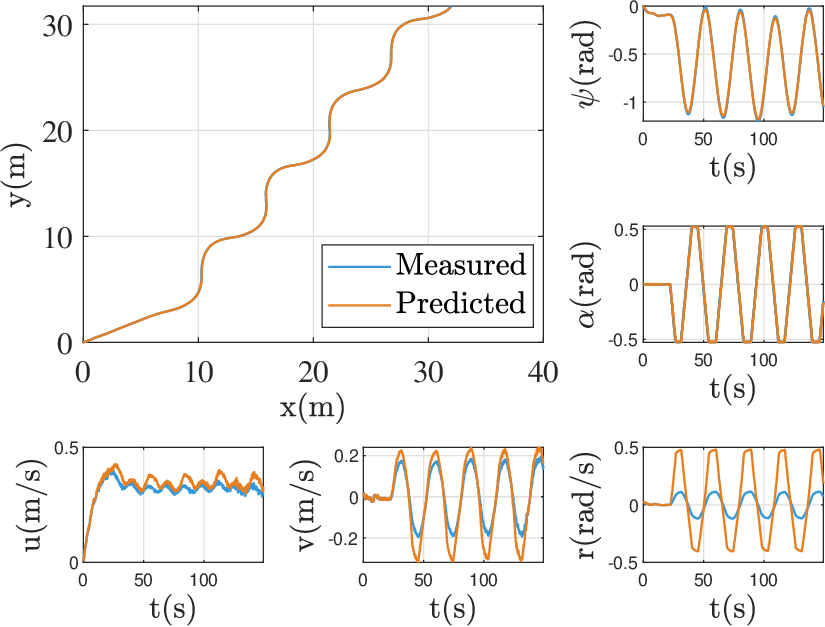}}\label{f RPS3_zigzag30_301}
 \hfill 	
  \subfloat[5 steps]{\includegraphics[width=0.3\textwidth]{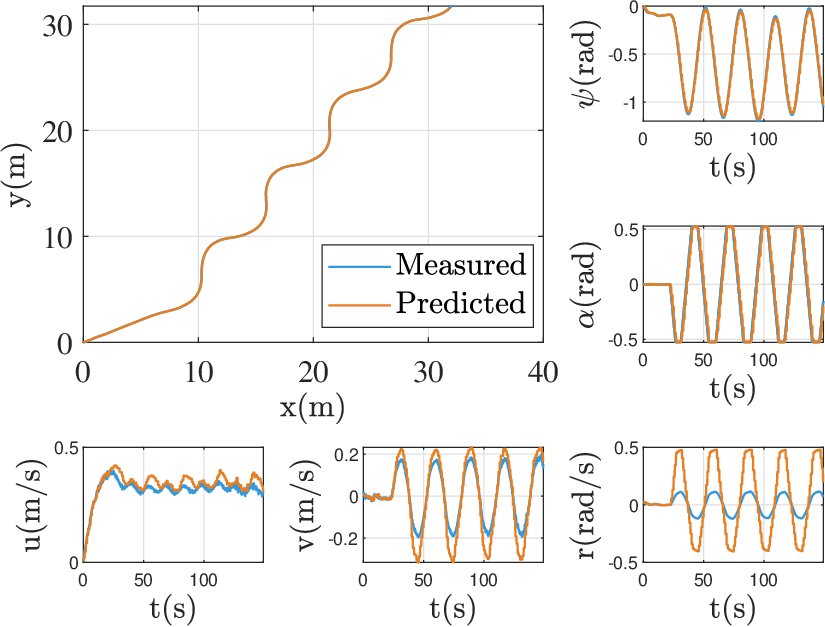}}\label{f RPS3_zigzag30_305}
 \hfill	
  \subfloat[15 steps]{\includegraphics[width=0.3\textwidth]{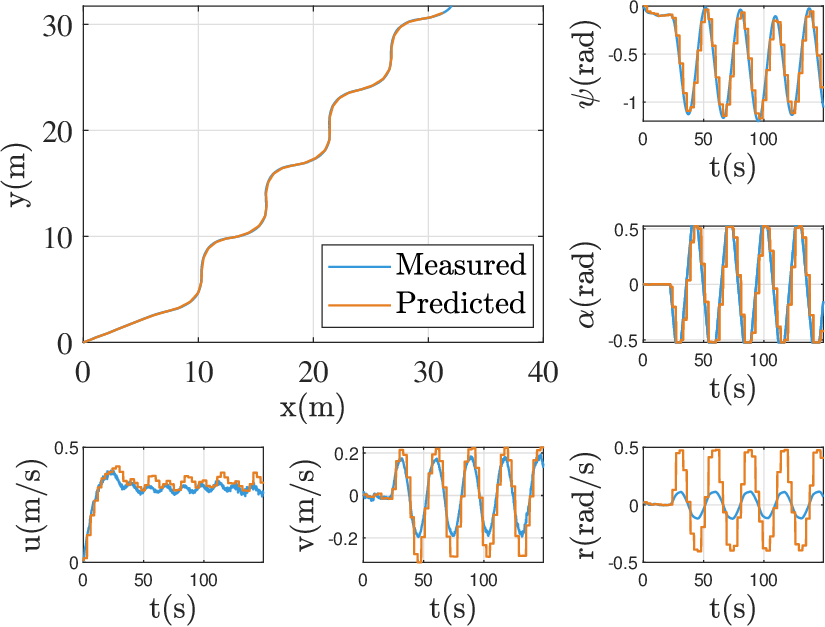}}\label{f RPS3_zigzag30_3015}
  \newline
  \subfloat[30 steps]{\includegraphics[width=0.3\textwidth]{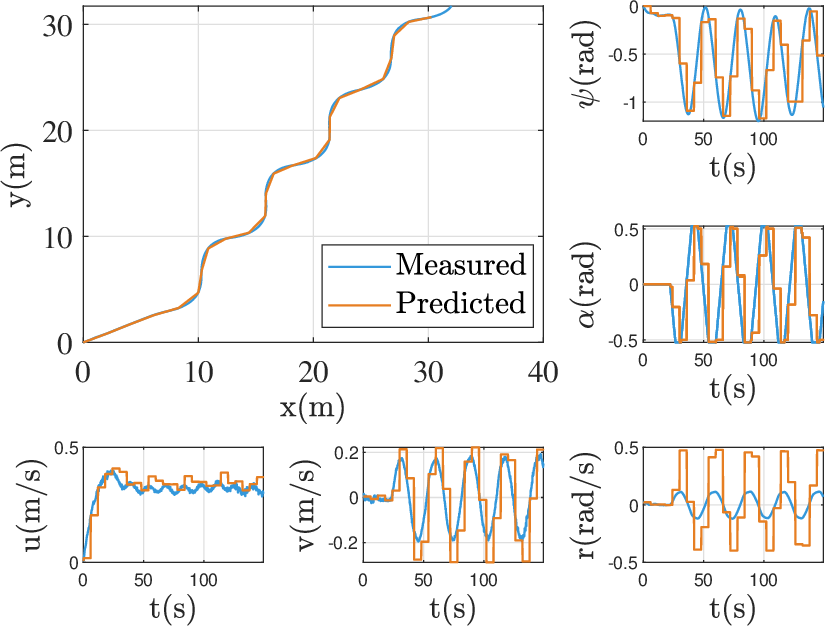}}
  \label{f RPS3_zigzag30_3030}
 \hfill 	
  \subfloat[50 steps]{\includegraphics[width=0.3\textwidth]{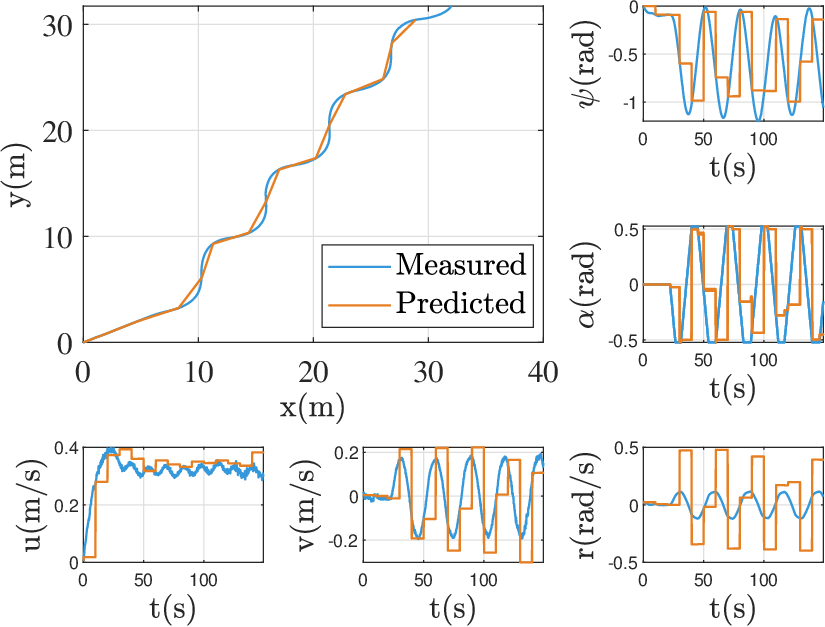}}
  \label{f RPS3_zigzag30_3050}
\caption{The validation based of different steps on +30+30 zigzag.}
\end{figure*}

\begin{figure*}[htbp]
  \subfloat[1 step]{\includegraphics[width=0.3\textwidth]{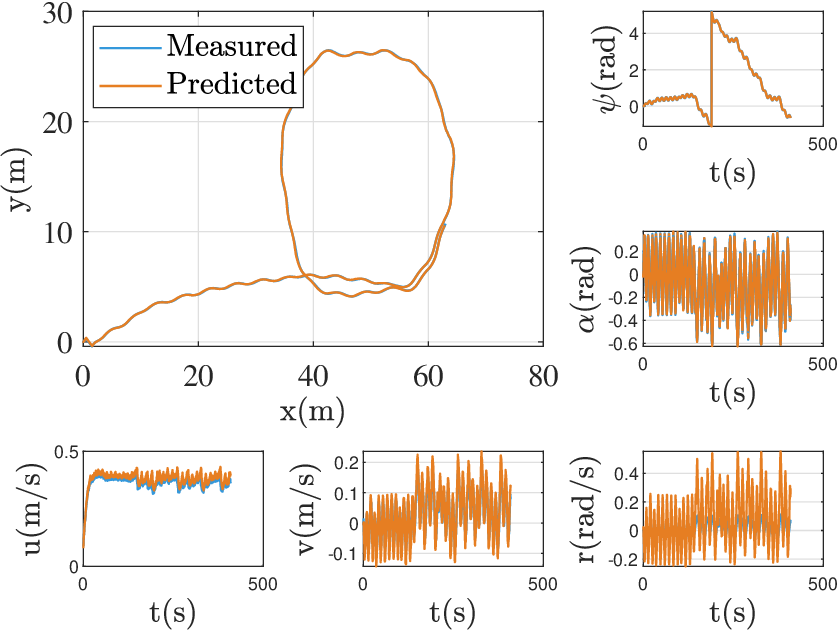}}\label{f RPS3_trajectroy_tracking1}
 \hfill 	
  \subfloat[5 steps]{\includegraphics[width=0.3\textwidth]{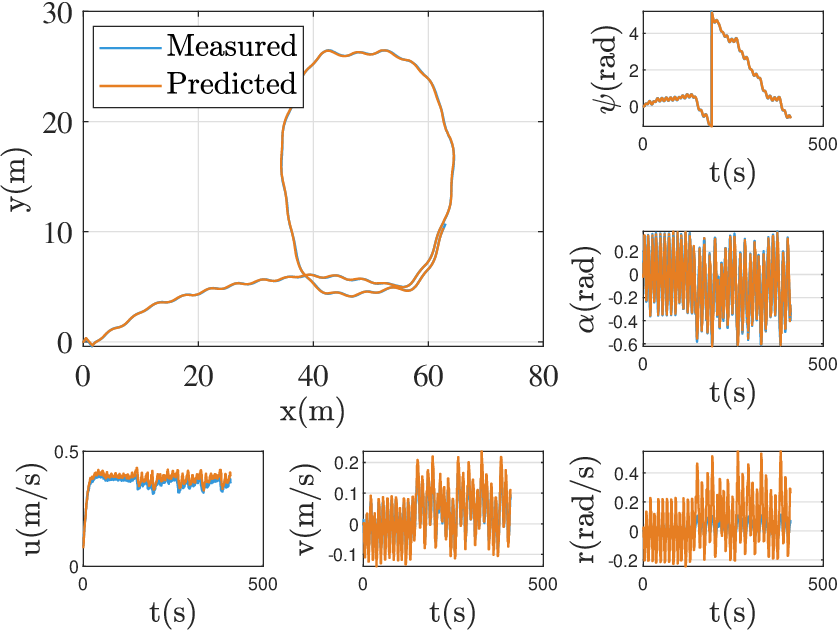}}\label{f RPS3_trajectroy_tracking5}
 \hfill	
  \subfloat[15 steps]{\includegraphics[width=0.3\textwidth]{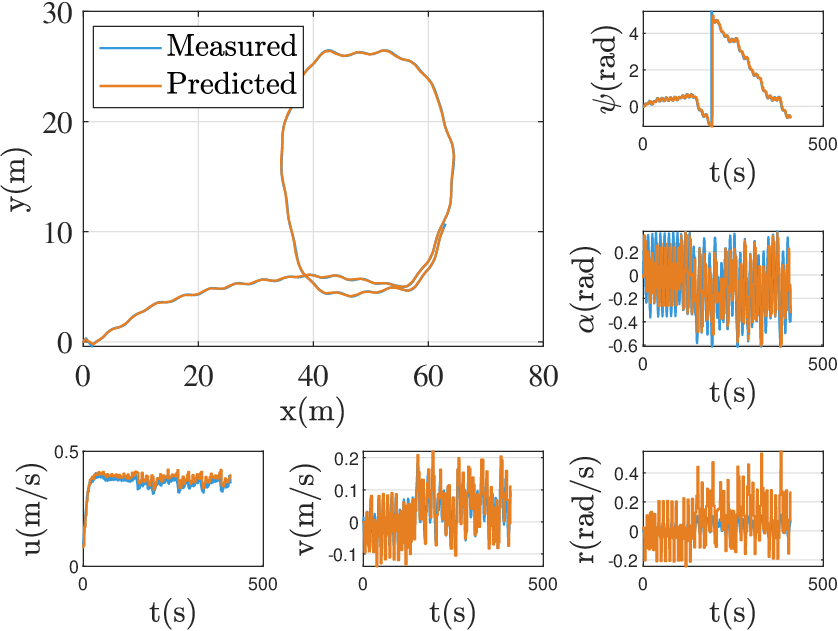}}\label{f RPS3_trajectroy_tracking15}
  \newline
  \subfloat[30 steps]{\includegraphics[width=0.3\textwidth]{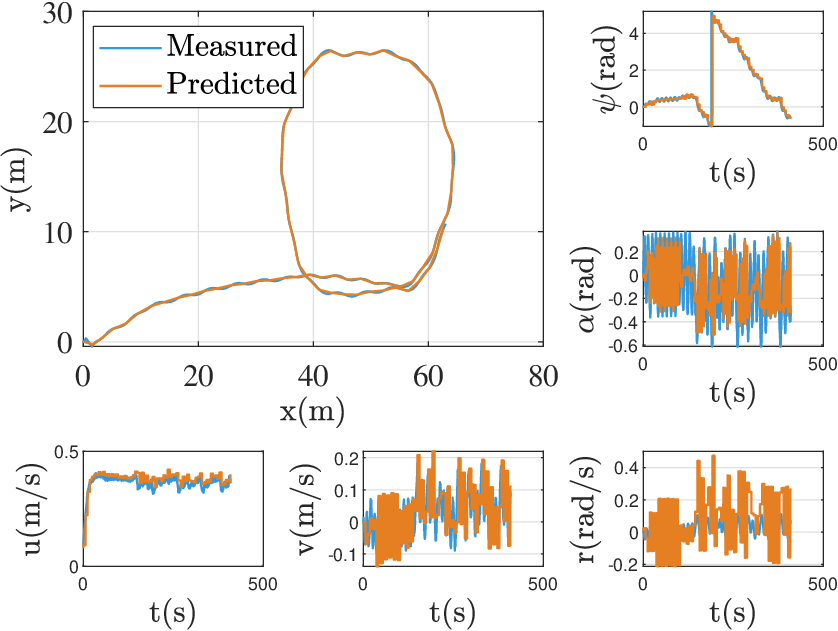}}
  \label{f RPS3_trajectroy_tracking30}
 \hfill 	
  \subfloat[50 steps]{\includegraphics[width=0.3\textwidth]{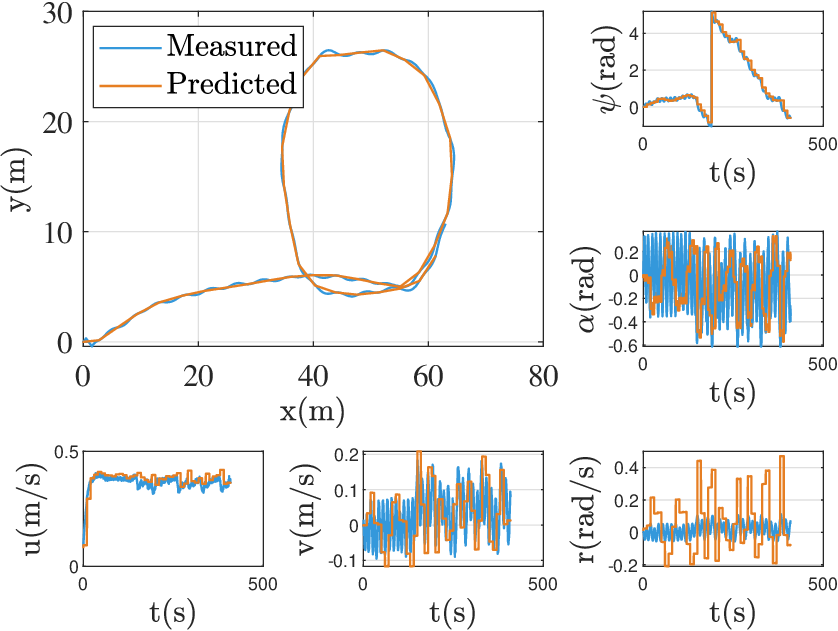}}
  \label{f RPS3_trajectroy_tracking50}
\caption{The validation based of different steps on +30+30 zigzag.}
\end{figure*}

In the previous subsection, we presented the validation for the ship's whole motion model with a single-step prediction. In this subsection, we will discuss the validation for multiple steps of prediction. We designed prediction cases for steps 1, 5, 10, 15, 20, 25, 30, 35, 40, 45, and 50 to evaluate the effectiveness of the proposed model parameters for the ship's movement. Figures \ref{f RPS3_zigzag30_301} to \ref{f RPS3_zigzag30_3050} show some of these prediction cases, where the input command is given by ourselves. Both cases exhibit the correct trajectory state trend; however, the error in the states increases as the step size of the given state increases, as demonstrated in Figure \ref{f RMSE of different steps}. According to the RMSE values of r, the predicted results are considered adequate for all the cases up to 40 steps. Even the 50-step prediction is acceptable for controller design as long as we can effectively compensate for the error. This demonstrates the effectiveness of the results obtained through GO, as they maintain a strong predictive performance even under different step sizes. Additionally, we validated the results using another trajectory state, as shown in Figures \ref{f RPS3_trajectroy_tracking1} to \ref{f RPS3_trajectroy_tracking50}. The RMSE results are consistent with the previous test, which further supports the effectiveness of the proposed method. These results are particularly useful in designing controllers, especially for model predictive control (MPC).

\section{Conclusion}
This study proposes a comprehensive solution for identifying the motion model of an under-actuated ship, from actual control input commands to ship motion states. We address the challenge of large-scale data in ship model identification by combining LO and GO methods and imposing constraints on ship features to limit the optimization solution domain. This approach allows us to utilize real data for large-scale nonlinear optimization solutions. Moreover, we use data from a real tug with two aft azimuth thrusters to estimate the ship's motion model. The proposed method successfully identifies the optimal solution for the target ship model parameters, and the predicted results align closely with the real measured values. These findings validate the effectiveness and sufficiency of the ship model construction and parameter estimation proposed in this study. 

In future work, we recommend considering the error in the yaw angle to improve yaw motion control. Additionally, we propose further exploration of model predictive control (MPC) as a suitable control method for this model to conduct more autonomous control experiments. Furthermore, we acknowledge that control algorithms tend to focus on stability and robustness, often neglecting the importance of accurate models. By building accurate models, the reliance on control algorithm design can be significantly reduced. Thus, investigating the balance between model accuracy and the effectiveness of control algorithms presents an interesting area for future research.


\printcredits

\bibliographystyle{cas-model2-names}

\bibliography{cas-refs}

\end{document}